\input ./style/arxiv-general.cfg
\documentclass[aop,MSNbibl,seceqn,dvips]{arximspdf}
\makeatletter
   \@ifpackageloaded{graphicx}{}{\usepackage{graphicx}}
\makeatother

\doi{10.1214/14-AOP914} 
\volume{43}
\issue{4}
\pubyear{2015}
\firstpage{1594}
\lastpage{1642}

\makeatletter
\newcommand{\rrVert}{\Vert}
\newcommand{\rrvert}{\vert}
\newcommand{\llVert}{\Vert}
\newcommand{\llvert}{\vert}
\newtheorem{teo}{Theorem}[section]
\newproclaim{defn}[teo]{Definition}
\newtheorem{cor}[teo]{Corollary}
\newtheorem{lem}[teo]{Lemma}
\newtheorem{propn}[teo]{Proposition}
\newproclaim{rem}[teo]{Remark}
\newproclaim{assump}[teo]{Assumption}
\makeatother

\begin{document}
\begin{frontmatter}

\title{Quenched invariance principles for random walks and elliptic
diffusions in random media with~boundary}
\runtitle{Quenched invariance principles for random media}

\begin{aug}
\author[A]{\fnms{Zhen-Qing} \snm{Chen}\thanksref{T1}\ead[label=e1]{zqchen@uw.edu}},
\author[B]{\fnms{David A.} \snm{Croydon}\ead[label=e2]{d.a.croydon@warwick.ac.uk}}
\and
\author[C]{\fnms{Takashi} \snm{Kumagai}\corref{}\thanksref{T2}\ead[label=e3]{kumagai@kurims.kyoto-u.ac.jp}}
\runauthor{Z.-Q. Chen, D. A. Croydon and T. Kumagai}
\affiliation{University of Washington, University of Warwick and Kyoto University}
\dedicated{Dedicated to Professors Martin T. Barlow and Ed. Perkins\\
on the occasion of their 60th birthdays}
\address[A]{Z.-Q. Chen\\
Department of Mathematics\\
University of Washington\\
Seattle, Washington 98195\\
USA\\
\printead{e1}} 
\address[B]{D. A. Croydon\\
Department of Statistics\\
University of Warwick\\
Coventry, CV4 7AL\\
United Kingdom\\
\printead{e2}}
\address[C]{T. Kumagai\\
Research Institute for Mathematical Sciences\\
Kyoto University\\
Kyoto 606-8502\\
Japan\\
\printead{e3}}
\end{aug}
\thankstext{T1}{Supported in part by NSF Grant DMS-12-06276
and NNSFC Grant 11128101.}
\thankstext{T2}{Supported in part by Grant-in-Aid for Scientific
Research (B) 22340017 and (A) 25247007.}

\received{\smonth{6} \syear{2013}}
\revised{\smonth{1} \syear{2014}}

%
\begin{abstract}
Via a Dirichlet form extension theorem and making full use of two-sided
heat kernel estimates, we establish quenched invariance principles for
random walks in random environments with a boundary. In particular, we
prove that the random walk on a supercritical percolation cluster or
among random conductances bounded uniformly from below in a half-space,
quarter-space, etc., converges when rescaled diffusively to a
reflecting Brownian motion, which has been one of the important open
problems in this area.
We establish a similar result for the random conductance model in a
box, which allows us to
improve existing asymptotic estimates for the relevant mixing time. Furthermore,
in the uniformly elliptic case, we present quenched invariance
principles for domains with more general boundaries.
\end{abstract}

%
\begin{keyword}[class=AMS]
\kwd[Primary ]{60K37}
\kwd{60F17}
\kwd[; secondary ]{31C25}
\kwd{35K08}
\kwd{82C41}
\end{keyword}
\begin{keyword}
\kwd{Quenched invariance principle}
\kwd{Dirichlet form}
\kwd{heat kernel}
\kwd{supercritical percolation}
\kwd{random conductance model}
\end{keyword}

\end{frontmatter}

\section{Introduction}\label{sec1}
Invariance principles for random walks in $d$-dimensional reversible
random environments date back to the 1980s \cite
{DFGW,KipVar,Kozlov,Kunn}. The most robust of the early results in this
area concerned
scaling limits for the annealed law, that is, the distribution of the
random walk averaged over the possible realizations of the environment,
or possibly established a slightly stronger statement involving some
form of convergence in probability. Studying the behavior of the random
walks under the quenched law, that is, for a fixed realization of the
environment, has proved to be a much more difficult task, especially
when there is some degeneracy in the model. This is because it is often
the case that a typical environment has ``bad'' regions that need to be
controlled. Nevertheless, over the last decade significant work has
been accomplished in this direction. Indeed, in the important case of
the random walk on the unique infinite cluster of supercritical (bond)
percolation on $\mathbb{Z}^d$, building on the detailed transition
density estimates of \cite{bar00}, a Brownian motion scaling limit has
now been established \cite{BerBis07,MatPia07,SidSzn04}. Additionally,
a number of extensions to more general random conductance models have
also been proved \cite{ABDH,BarDeu08,BiskP,Mat}.

%
\begin{figure}[b]

\includegraphics{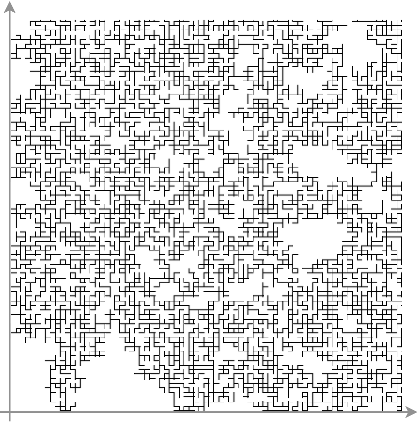}

\caption{A section of the unique infinite cluster for supercritical
percolation on $\mathbb{Z}_+^2$ with parameter $p=0.52$.}\label{intfig}
\end{figure}

While the above body of work provides some powerful techniques for
overcoming the technical challenges involved in proving quenched
invariance principles, such as studying ``the environment viewed from
the particle'' or the ``harmonic corrector'' for the walk (see \cite
{Biskup} for a survey of the recent developments in the area), these
are not without their limitations. Most notably, at some point, the
arguments applied all depend in a fundamental way on the translation
invariance or ergodicity under random walk transitions of the
environment. As a consequence, some natural variations of the problem
are not covered. Consider, for example, supercritical percolation in a
half-space $\mathbb{Z}_+\times\mathbb{Z}^{d-1}$, or possibly an
orthant of $\mathbb{Z}^d$. Again, there is a unique infinite cluster
(Figure~\ref{intfig} shows a simulation of such in the first quadrant
of $\mathbb{Z}^2$), upon which one can define a random walk. Given the
invariance principle for percolation on $\mathbb{Z}^d$, one would
reasonably expect that this process would converge, when rescaled
diffusively, to a Brownian motion reflected at the boundary. After all,
as is illustrated in the figure, the ``holes'' in the percolation
cluster that are in contact with the boundary are only on the same
scale as those away from it. However, the presence of a boundary means
that the translation invariance/ergodicity properties necessary for
applying the existing arguments are lacking. For this reason, it has
been one of the important open problems in this area to prove the
quenched invariance principle for random walk on a percolation cluster,
or among random conductances more generally, in a half-space (see \cite
{BerBis07}, Section~B and \cite{Biskup}, Problem 1.9). Our aim is to
provide a new approach for overcoming this issue, and thereby establish
invariance principles within some general framework that includes
examples such as those just described.

The approach of this paper is inspired by that of \cite{BC1,BC2},
where the invariance principle for random walk on grids inside a given
Euclidean domain $D$ is studied. It is shown first in \cite{BC1} for a
class of bounded domains including Lipschitz domains and then in \cite
{BC2} for \textit{any} bounded domain $D$ that the simple random walk
converges to the (normally)
reflecting Brownian motion on $D$ when the mesh size
of the grid tends to zero. Heuristically, the normally
reflecting Brownian motion is a continuous Markov process on $\overline
{D}$ that behaves like Brownian motion in $D$ and is ``pushed back''
instantaneously along the inward normal direction when it hits the
boundary. See Section~\ref{sect2} for
a precise definition and more details. The main idea and approach of
\cite{BC1,BC2} is as follows:
(i) show that random walk killed upon hitting the boundary converges
weakly to the absorbing Brownian motion in $D$, which is trivial;
(ii) establish tightness for the law of random walks;
(iii) show any sequential limit is a symmetric Markov process and can
be identified with reflecting Brownian motion via a Dirichlet form
characterization. In \cite{BC1,BC2}, (ii) is achieved by using a
forward--backward martingale decomposition of the process and
the identification in (iii) is accomplished by using a result from
the boundary theory of Dirichlet form, which says that the reflecting
Brownian motion on $D$ is the maximal Silverstein's extension of the
absorbing Brownian motion in $D$; see~\cite{BC1}, Theorem~1.1 and
\cite{CF}, Theorem 6.6.9.

For quenched invariance principles for random walks in random
environments with a boundary, step (i) above can be established by
applying a quenched invariance principle for the full-space case. For
step (ii), that is, establishing tightness, the forward--backward
martingale decomposition method does not work well with unbounded
random conductances.
To overcome this difficulty,
as well as for the desire to establish an invariance principle for
every starting point, we will make the full use of detailed two-sided
heat kernel estimates for random walk on random clusters. In
particular, we provide sufficient conditions for the subsequential
convergence that involve the H\"{o}lder continuity of harmonic
functions (see Section~\ref{suffsec}). This continuity property can be
verified in examples
by using existing two-sided heat kernel bounds. We remark that the
corrector-type methods for full-space models, such as the approach of
\cite{BerBis07}, often require only upper bounds on the heat kernel.
Using the H\"older regularity, we can further show that any
subsequential limit of random walks in random environments is a
conservative symmetric Hunt process with continuous sample paths. In
step (iii), we can
identify the subsequential limit process with the reflecting Brownian motion
by a Dirichlet form argument (see Theorem~\ref{teosilvs}).
In summary,
our approach for proving quenched invariance principles for random
walks in random environments with a boundary encompasses two novel
aspects: a Dirichlet form extension argument and the full use of
detailed heat kernel estimates.

The full generality of the random conductance model to which we are
able to apply the above argument is presented in Section~\ref{percsec}.
As an illustrative application of Theorem~\ref{teosilvs}, though we
state here a theorem that verifies the conjecture described above
concerning the diffusive behavior of the random walk on a supercritical
percolation cluster on a half-space, quarter-space, etc. We recall that
the variable speed random walk (VSRW) on a connected (unweighted) graph
is the continuous time Markov process that jumps from a vertex at a
rate equal to its degree to a uniformly chosen neighbor (see
Section~\ref{percsec} for further details). In this setting, similar
results to that stated can be obtained for the so-called constant speed
random walk (CSRW), which has mean one exponential holding times, or
the discrete time random walk (see Remark~\ref{genrem} below).

Let $\mathbb{Z}_+:=\{0,1,2,\ldots\}$ and $\mathbb{R}_+:=[0,\infty)$.
Then the following is our main theorem.

\begin{teo} Fix $d_1,d_2\in\mathbb{Z}_+$ such that $d_1\geq1$ and
$d:=d_1+d_2\geq2$. Let $\mathcal{C}_1$ be the unique infinite cluster
of a supercritical bond percolation process on \mbox{$\mathbb
{Z}^{d_1}_+\times
\mathbb{Z}^{d_2}$}, and let $Y=(Y_t)_{t\geq0}$ be the associated VSRW.
For almost-every realization of $\mathcal{C}_1$, it holds that the
rescaled process $Y^n=(Y^n_t)_{t\geq0}$, as defined by
\[
Y^n_t:=n^{-1}Y_{n^2t},
\]
started from $Y^n_0=x_n\in n^{-1}\mathcal{C}_1$, where $x_n\rightarrow
x\in\mathbb{R}^{d_1}_+\times\mathbb{R}^{d_2}$, converges in
distribution to $\{X_{ct}; t\geq0\}$, where $c\in(0,\infty)$ is a
deterministic constant
and $\{X_t; t\geq0\}$ is the (normally) reflecting Brownian motion on
$\mathbb{R}^{d_1}_+\times\mathbb{R}^{d_2}$ started from $x$.
\end{teo}

As an alternative to unbounded domains, one could consider compact
limiting sets, replacing $Y^n$ in the previous theorem by the rescaled
version of the variable speed random walk on the largest percolation
cluster contained in a box $[-n,n]^d\cap\mathbb{Z}^d$, for example. As
presented in Section~\ref{boxsec}, another application of Theorem~\ref{teosilvs} allows an invariance principle to be established in this
case as well, with the limiting process being Brownian motion in the
box $[-1,1]^d$, reflected at the boundary. Consequently, we are able to
refine the existing knowledge of the mixing time asymptotics for the
sequence of random graphs in question from a tightness result \cite
{BenMos} to an almost-sure convergence one (see Corollary~\ref{mtcor} below).

Although in the percolation setting we only consider relatively simple
domains with ``flat'' boundaries, this is mainly for technical reasons
so that deriving the percolation estimates in Section~\ref{percosect}
required for our proofs is manageable. Indeed, in the case when we
restrict to uniformly elliptic random conductances, so that controlling
the clusters of extreme conductances is no longer an issue, we are able
to derive from Theorem~\ref{teosilvs} quenched invariance principles
in any uniform domain, the collection of which forms a large class of
possibly nonsmooth domains that includes (global) Lipschitz domains
and the classical van Koch snowflake planar domain as special cases.
These applications are discussed in Section~\ref{uesec}.

Homogenization of reflected SDE/PDE on half-planes and more general
domains has been studied in various contexts (see, e.g., \cite
{BLLS,BLP,KLO,Rho,Tana}; we refer to \cite{JKO,KLO} and the references
therein for the history of homogenization for diffusions in random
environments). In a recent paper \cite{Rho}, Rhodes proves
homogenization (as a convergence in product measure
in environment and state space of quenched distribution,
which implies an annealed invariance principle) for symmetric reflected
diffusions in upper half-spaces.
His method is based on the Girsanov formula and a use of subsidiary
diffusions with an invariant probability measure, which is very
different from ours. Although we can also only handle symmetric cases,
our methods contribute to this field as well. This is because the
analytical part of our results (namely Section~\ref{sect2}) holds
for the entire class of uniform domains. Moreover, our results are on
the level of quenched invariance principles. The presentation of how
our techniques can be applied in the uniformly elliptic random
divergence form setting appears in Section~\ref{dfsec}. Note further
that in this setting we resolve the open problem on the quenched
invariance principle starting from arbitrary starting points posed in
\cite{Rho}, pages 1004--1005.

The remainder of the paper is organized as follows. In Section~\ref
{sect2}, we introduce an abstract framework for proving invariance
principles for reversible Markov processes in a Euclidean domain. This
is applied in Section~\ref{percsec} to our main example of a random
conductance model in half-spaces, quarter-spaces, etc. The details of
the other examples discussed above are presented in Section~\ref
{exsec}. Our results for the random conductance model depend on a
number of technical percolation estimates, some of the proofs of which
are contained in the \hyperref[sec5]{Appendix} that appears at the end of this article.
The \hyperref[sec5]{Appendix} also contains a proof of a generalization of existing
quenched invariance principles that allows for arbitrary starting
points (previous results have always started the relevant processes
from the origin, which will not be enough for our purposes).

Finally, in this paper, for a locally compact separable metric space $E$,
we use $C_b(E)$ and $C_\infty(E)$ to denote the space
of bounded continuous functions on $E$ and the space of continuous
functions on $E$ that vanish at infinity, respectively.
The space of continuous functions on $E$ with compact support
will be denoted by $C_c(E)$. For real numbers $a, b$, we use
$a\vee b$ and $a\wedge b$ for $\max\{a, b\}$ and $\min
\{a, b\}$, respectively.

\section{Framework}\label{sect2}

The following definition is taken from V\"ais\"al\"a \cite{Va},
where various equivalent definitions are discussed.
An open connected subset $D$ of $\mathbb{R}^d$ is called
\textit{uniform} if there exists a constant $C$
such that for every $x,y\in D$ there is a rectifiable curve $\gamma$
joining $x$ and $y$ in $D$ with
length $(\gamma)\leq C|x-y|$ and moreover
$\min\{ |x-z|, |z-y|
\}\leq C \operatorname{dist} (z, \partial D)$
for all points $z\in\gamma$.
Here, $\operatorname{dist} (z, \partial D)$ is the Euclidean distance
between the point $z$ and the set $\partial D$.
Note that a uniform domain
with respect to an inner metric
is called \textit{inner uniform} in \cite{GSC}, Definition 3.6.

For example, the classical van Koch snowflake
domain in the conformal mapping theory is a uniform domain
in $\mathbb{R}^2$. Every (global) Lipschitz domain is uniform,
and every
\textit{nontangentially accessible domain} defined by Jerison and
Kenig in \cite{JK} is a uniform domain (see (3.4) of \cite{JK}).
However, the boundary of a uniform domain can be highly
nonrectifiable and, in general, no regularity of its boundary
can be inferred (besides the easy fact that the Hausdorff dimension
of the boundary is strictly less than $d$).

It is known (see Example 4 on page 30
and Proposition~1 in Chapter VIII of~\cite{JW})
that any uniform domain in $\mathbb{R}^d$ has $m(\partial D)=0$
and there exists a positive constant $c>0$
such that
%
\begin{equation}
\label{e21n} m \bigl(D \cap B_E(x, r)\bigr)\geq c r^n
\qquad\mbox{for all } x\in\overline{D}\mbox{ and }0<r\leq1,
\end{equation}
where $m$ denotes the Lebesgue measure in
$\mathbb{R}^d$ and $B_E(x,r)$ denotes the Euclidean ball of radius $r$
centered at $x$.

Let $D$ be a uniform domain in $\mathbb{R}^d$.
Suppose $(A(x))_{x\in\overline{D}}$ is a measurable symmetric
$d\times
d$ matrix-valued function such that
%
\begin{equation}
\label{e22n} c^{-1}I\leq A(x) \leq cI \qquad\mbox{for a.e. } x\in
\overline{D},
\end{equation}
where $I$ is the $d$-dimensional identity matrix and $c$ is a constant
in $[1,\infty)$. Let
%
\begin{equation}
\label{e23} \mathcal{E}(f,g):=\frac{1}{2} \int_D
\nabla f(x) \cdot A(x)\nabla g(x)\,dx \qquad\mbox{for } f, g\in W^{1,2}(D),
\end{equation}
where
\[
W^{1,2}(D):= \bigl\{f\in L^2(D; m)\dvtx \nabla f\in
L^2(D; m) \bigr\}. %
\]
An important property of a uniform domain $D\subset\mathbb{R}^d$
is that there is a bounded linear extension operator
$T\dvtx W^{1,2}(D) \to W^{1,2}(\mathbb{R}^d)$ such that
$Tf=f$ a.e. on $D$ for $f\in W^{1,2}(D)$.
It follows that $(\mathcal{E}, W^{1,2}(D))$ is
a regular Dirichlet form on $L^2(\overline{D}; m)$ and so
there is a continuous diffusion process $X=(X_t, t\geq0; \mathbb
{P}_x,x\in
\overline{D})$
associated with it, starting from $\mathcal{E}$-quasi-every point.
Here, a property is said to hold $\mathcal{E}$-quasi-everywhere
means that there is a set ${\mathcal N}\subset\overline{D}$
having zero capacity with respect to
the Dirichlet form $(\mathcal{E}, W^{1,2}(D))$ so that the property holds
for points in~${\mathcal N}^c$.
According to \cite{GSC}, Theorem 3.10 and (\ref{e21n})
(see also \cite{BCR}, (3.6)),
$X$ admits a jointly continuous transition density function
$p(t, x, y)$ on $\mathbb{R}_+ \times\overline{D} \times\overline{D}$
and
%
\begin{equation}
\label{e23n} \qquad {c_1} {t^{-d/2}} \exp\biggl(-\frac{c_2|x- y|^2}{t}
\biggr) \leq p(t, x, y)\leq{c_3} {t^{-d/2}} \exp\biggl(-
\frac{c_4 |x-y|^2}{t} \biggr)
\end{equation}
for every $x, y\in\overline{D}$ and $0<t\leq1$.
Here, the constants $c_1,\ldots, c_4>0$
depend on the diffusion matrix $A(x)$ only through the ellipticity
bound $c$ in (\ref{e22n}).
Consequently, $X$ can be refined so that it can start from every point
in $\overline{D}$.
The process $X$ is called a symmetric reflecting diffusion on
$\overline{D}$.
We refer to \cite{C} for sample path properties\vspace*{1pt} of $X$. When $A=I$, $X$
is the (normally) reflecting Brownian motion on $\overline{D}$.
Reflecting Brownian motion $X$ on $\overline{D}$ in general
does not need to be semi-martingale.
When $\partial D$ locally has finite lower Minkowski content, which is
the case when
$D$ is a Lipschitz domain, $X$ is a semi-martingale and admits
the following Skorohod decomposition
(see \cite{C2}, Theorem 2.6):
%
\begin{equation}
\label{e25} X_t =X_0 +W_t+\int
_0^t \vec n (X_s)
\,dL_s,\qquad t\geq0.
\end{equation}
Here, $W$ is the standard Brownian motion in $\mathbb{R}^d$,
$\vec n$ is the unit inward normal vector field of $D$ on $\partial D$,
and $L$ is a positive continuous additive functional of $X$ that
increases only when $X$ is on the boundary, that is, $L_t=\int_0^t
1_{\{
X_s\in\partial D\}} \,dL_s$ for $t\geq0$.
Moreover, it is known that the reflecting Brownian motion spends zero
Lebesgue amount of time at
the boundary $\partial D$. These together with (\ref{e25}) justify
the heuristic description we gave in the \hyperref[sec1]{Introduction} for the
reflecting Brownian motion in $D$.

\subsection{Convergence to reflecting diffusion}\label{sec2.1}

In this subsection, $D$ is a uniform domain in $\mathbb{R}^d$ and $X$
is a
reflecting diffusion process on $\overline{D}$ associated
with the Dirichlet form $(\mathcal{E}, W^{1, 2}(D))$ on $L^2(D; m)$
given by
(\ref{e23}).
Denote by \mbox{$(X^D,\mathbb{P}_x^D,x\in\overline{D})$}
the subprocess of $X$ killed on exiting $D$. It is known (see, e.g.,
\cite{CF})
that the Dirichlet form of $X^D$ on $L^2(D; m)$ is $(\mathcal{E},
W^{1,2}_0(D))$, where
\[
W^{1,2}_0(D):= \bigl\{f\in W^{1,2}(D)\dvtx f=0\
\mathcal{E}\mbox{-quasi-everywhere on } \partial D \bigr\}.
\]

Suppose that $\{D_n; n\ge1\}$ is a sequence of Borel subsets of
$\overline{D}$ such that
each $D_n$ supports a measure $m_n$ that
converges vaguely to the Lebesgue measure $m$
on $\overline{D}$.
The following result plays a key role in our approach to the
quenched invariance principle for random walks in random environments
with boundary.

\begin{teo}\label{teosilvs} For each $n\in\mathbb N$, let $(X^n,\mathbb{P}
^{n}_{x},x\in D_n)$ be an $m_n$-symmetric Hunt process on $D_n$.
Assume that for every subsequence $\{n_j\}$, there exists a
sub-subsequence $\{n_{j(k)}\}$ and a continuous
conservative $m$-symmetric strong Markov
process $(\widetilde X, \widetilde{\mathbb{P}}_x, x\in\overline{D})$
such that the following
three conditions are satisfied:
\begin{longlist}[(iii)]
\item[(i)] for every $x_{n_{j(k)}}\to x$ with $x_{n_{j(k)}}\in
D_{n_{j(k)}}$, $\mathbb{P}^{n_{j(k)}}_{x_{n_{j(k)}}}$
converges weakly in $\mathbb{D}([0,\infty),\overline{D})$ to
$\widetilde{\mathbb{P}}_x$;\vspace*{1pt}
\item[(ii)] $\widetilde X^D$, the subprocess of $\widetilde X$
killed upon
leaving $D$, has the same distribution as $X^D$;
\item[(iii)] the Dirichlet form $(\widetilde{\mathcal{E}}, \widetilde
{\mathcal{F}})$
of $\widetilde{X}$ on $L^2(D; m)$ has the properties that
%
\begin{equation}
\label{e21a} {\mathcal C}\subset\widetilde{\mathcal{F}}\quad\mbox{and}\quad\widetilde{\mathcal{E}}(f,f)\le C_0 \mathcal{E}(f,f)\qquad
\mbox{for every } f\in{
\mathcal C},
\end{equation}
where $\mathcal{C}$ is a core for the Dirichlet form $(\mathcal{E},
W^{1,2}(D))$ and
$C_0\in[1,\infty)$ is a constant.
\end{longlist}
It then holds that for every $x_n\rightarrow x$ with $x_n\in D_n$,
$(X^{n}, \mathbb{P}^n_{x_n})$ converges weakly in $\mathbb
{D}([0,\infty),\overline{D})$ to $(X, \mathbb{P}_x)$.
\end{teo}

\begin{pf}
With both $\widetilde X$ and $X$ being
$m$-symmetric Hunt processes
on $\overline{D}$, it suffices to show that their corresponding\vspace*{1pt} (quasi-regular)
Dirichlet forms on $L^2 (D; m)$ are the same; that is
$(\widetilde{\mathcal{E}}, \widetilde{\mathcal{F}})=(\mathcal{E},
W^{1,2}(D))$. Condition (iii) immediately
implies that $W^{1,2}(D)\subset\widetilde{\mathcal{F}}$ and
\[
\widetilde{\mathcal{E}}( f, f) \leq C_0 \mathcal{E}(f, f) \qquad
\mbox{for every } f\in W^{1,2}(D).
\]
Next, observe that since $\widetilde X$ is a\vspace*{2pt} diffusion process
admitting no
killings, its associated Dirichlet form is strongly local. Thus,\vspace*{1pt} for
every $u\in\widetilde{\mathcal{F}}$, $\widetilde{\mathcal{E}}(u,
u)=\frac{1}2 \tilde\mu_{\langle u\rangle
}(\overline{D})$,
where $\tilde\mu_{\langle u\rangle}$ is the energy measure
corresponding to $u$. By
the proof of \cite{Mos}, proposition on page 389,
%
\begin{eqnarray}
\tilde\mu_{\langle u\rangle}(dx) \leq C_0\nabla u(x) A(x)\nabla
u(x)\,dx\leq cC_0 \bigl|\nabla u (x)\bigr|^2 \,dx
\nonumber
\\
\eqntext{\mbox{on }\overline{D} \mbox{ for every } u\in W^{1,2}(D).}
\end{eqnarray}
This in particular implies that
%
\begin{equation}
\label{e21} \tilde\mu_{\langle u\rangle}(\partial D) =0 \qquad\mbox
{for } u
\in W^{1,2}(D).
\end{equation}
On the other hand, by the strong local property of $\tilde\mu
_{\langle u\rangle}$
and the fact that $\widetilde X^D$ has the same distribution as $X^D$,
we have
that every bounded function in $\widetilde{\mathcal{F}}$---the
collection of which we
denote by $\widetilde{\mathcal{F}}_b$---is locally in $ W^{1,2}_0(D)$ and
%
\begin{equation}
\label{e22} \mathbf{1}_D (x) \tilde\mu_{\langle u\rangle}(dx) =
\mathbf{1}_D (x) \nabla u(x) A(x)\nabla u(x)\,dx \qquad\mbox{for } u\in
\widetilde{\mathcal{F}}_b.
\end{equation}
This together with (\ref{e21}) implies that $\widetilde{\mathcal{E}}(u,
u)=\mathcal{E}(u,
u)$ for every bounded $u\in W^{1,2}(D)$, and\vspace*{1pt} hence for every $u\in
W^{1,2}(D)$. Furthermore, (\ref{e22}) implies that for $u\in
\widetilde{\mathcal{F}}
_b$, $\int_D |\nabla u (x)|^2 \,dx <\infty$ and so $u\in W^{1,2}(D)$.
Consequently, we have $\widetilde{\mathcal{F}}\subset W^{1,2}(D)$, and
thus $(\widetilde{\mathcal{E}},
\widetilde{\mathcal{F}})=(\mathcal{E}, W^{1,2}(D))$.
\end{pf}

\begin{rem}\label{teoconserrem}
(i) Note that if $(X^n_t)_{t\ge0}$ is conservative for each $n\in
\mathbb{N}$ and
$\{\mathbb{P}^{n_{j(k)}}_{x_{n_{j(k)}}}\}$ is tight, then $\widetilde
X$ is
conservative.\vspace*{1pt}

(ii) Theorem~\ref{teosilvs} can be viewed as a variation of \cite{BC1},
Theorem 1.1.
The difference is that in \cite{BC1}, Theorem 1.1, the constant $C_0$
in (\ref{e21a}) is assumed to be 1 but the limiting process
$\widetilde X$
only need to be Markov and does not need to be continuous a priori,
while for Theorem~\ref{teosilvs}, the condition on the constant $C_0$
is weaker but we need to assume a priori that the limit process
$\widetilde X$ is continuous.
\end{rem}

\subsection{Sufficient condition for subsequential convergence}\label{sec2.2}\label{suffsec}

In this subsection, we give some sufficient conditions for the
subsequential convergence of $\{X^n\}$; in other words, sufficient
conditions for (i) in Theorem~\ref{teosilvs}. For simplicity, we
assume that $0\in D_n$ for all $n\ge1$ throughout this section, though
note this restriction can easily be removed.

We start by introducing our first main assumption, which will allow us
to check an equi-continuity property for the $\lambda$-potentials
associated with the elements of~$\{X^n\}$ (see Proposition~\ref{PR3}
below). In the statement of the assumption, we suppose that $(\delta
_n)_{n\geq1}$ is a decreasing sequence in $[0,1]$ with $\lim_{n\to
\infty}\delta_n=0$ and such that $|x-y|\ge\delta_n$ for all distinct
$x,y\in D_n$. (When $\delta_n\equiv0$, this condition always holds.
However, our assumption will give an additional restriction.) We denote by~$\tau_A(X^n)$ the first exit time of the process $X^n$ from the set $A$.

\begin{assump}\label{teoassumpHC} There\vspace*{1pt} exist $c_1, c_2, c_3, \beta,\gamma\in
(0,\infty)$, $N_0\in\mathbb{N}$ such that the following
hold for all $n\ge N_0$, $x_0\in B_E(0,c_1n^{1/2})$,
and $\delta_n^{1/2}\le r\le1$.
\begin{longlist}[(ii)]
\item[(i)] For all $x\in B_E(x_0,r/2)\cap D_n$,
\[
\mathbb{E}^n_x \bigl[ \tau_{B_E(x_0,r)\cap D_n}
\bigl(X^n\bigr) \bigr]\leq c_2r^{\beta}.
\]

\item[(ii)] If $h_n$ is bounded in $D_n$ and harmonic (with respect
to $X^n$) in
a ball $B_E(x_0,r)$, then
\[
\bigl|h_n(x)-h_n(y)\bigr|\leq c_3 \biggl(
\frac{|x-y|}{r} \biggr)^\gamma\|h_n\| _\infty
\qquad\mbox{for } x,y\in B_E(x_0,r/2)\cap
D_n.
\]
\end{longlist}
\end{assump}

Define for $\lambda>0$ the $\lambda$-potential
\[
U_n^\lambda f(x)=\mathbb{E}^n_x\int
_0^\infty e^{-\lambda t} f\bigl(X^n_t
\bigr) \,dt \qquad\mbox{for } x\in D_n.
\]

\begin{propn}\label{PR3} Under Assumption~\ref{teoassumpHC}, there exist
$C=C_\lambda\in(0,\infty)$ and $\gamma'\in(0,\infty)$ such that the
following holds for any bounded function $f$ on $D_n$,
for any $n\geq N_0$ and any $x,y\in D_n$ such that $x\in
B_E(0,c_1n^{1/2})$ and $|x-y|< 1/4$:
%
\begin{equation}
\label{cbe000} \bigl|U_n^\lambda f(x)-U_n^\lambda
f(y)\bigr|\leq C|x-y|^{\gamma'}\|f\|_\infty.
\end{equation}
In particular, we have
%
\begin{equation}
\label{cbe00110} \lim_{\delta\to0}\ \sup_{n\geq N_0}\ \mathop{\sup_{x,y\in
D_n\cap
B_E(0,c_1n^{1/2})\dvtx}}_{|x-y|<\delta}\bigl|U_n^\lambda
f(x)-U_n^\lambda f(y)\bigr|=0.
\end{equation}
\end{propn}

\begin{pf}
The proof is similar to that of \cite{BKK}, Proposition 3.3.
Fix $x_0\in B_E(0,c_1n^{1/2})\cap D_n$, let $1\ge r\ge\delta_n^{1/2}$,
and suppose $x,y\in B_E(x_0,r/2)$.
Set $\tau^n_r:=\tau_{B_E(x_0,r)\cap D_n}(X^n)$. By the strong Markov property,
\begin{eqnarray*}
U_n^\lambda f(x)&=&\mathbb{E}^n_x
\int_0^{\tau^n_r} e^{-\lambda t} f
\bigl(X_t^n\bigr) \,dt +\mathbb{E}^n_x
\bigl[ \bigl(e^{-\lambda\tau^n_r}-1\bigr)U_n^\lambda f
\bigl(X^n_{\tau^n_r}\bigr) \bigr] + \mathbb{E}^n_x
\bigl[ U_n^\lambda f\bigl(X^n_{\tau^n_r}
\bigr) \bigr]
\\
&=&I_1+I_2+I_3
\end{eqnarray*}
and similarly when $x$ is replaced by $y$. We have by Assumption~\ref{teoassumpHC}(i) that
\[
|I_1|\leq\|f\|_\infty\mathbb{E}^n_x
\tau^n_r\leq c_2r^{\beta}\|f\|
_\infty
\]
and by noting $\|U_n^\lambda f\|_\infty\leq\frac{1}{\lambda}\|f\|
_\infty$ that
\[
|I_2|\leq\lambda\mathbb{E}^n_x
\tau^n_r \bigl\|U_n^\lambda f
\bigr\|_\infty\leq c_2r^{\beta}\| f\|_\infty.
\]
Similar statements also hold when $x$ is replaced by $y$. So,
%
\begin{equation}
\label{C05} \qquad\quad \bigl\llvert U_n^\lambda f(x)-U_n^\lambda
f(y)\bigr\rrvert\leq4c_2r^{\beta}\| f\|_\infty+
\bigl\llvert\mathbb{E}^n_x U_n^\lambda
f\bigl(X^n_{\tau^n_r}\bigr) -\mathbb{E}^n_y
U_n^\lambda f\bigl(X^n_{\tau^n_r}\bigr)\bigr
\rrvert.
\end{equation}
But $z\to\mathbb{E}^n_z U_n^\lambda f(X^n_{\tau^n_r})$ is bounded
in $\mathbb{R}^d$ and
harmonic in $B_E(x_0,r)$, so by Assumption~\ref{teoassumpHC}(ii), the
second term in (\ref{C05}) is bounded by $c_3(|x-y|/r)^\gamma\|
U_n^\lambda f\|_\infty$.
So by $\|U_n^\lambda f\|_\infty\leq\frac{1}{\lambda}\|f\|_\infty$
again, we have
%
\begin{eqnarray}\label{cbe01}
\bigl\llvert U_n^\lambda f(x)-U_n^\lambda
f(y)\bigr\rrvert\leq c \biggl(r^{\beta
}+\lambda^{-1} \biggl(
\frac{|x-y|}r \biggr)^{\gamma} \biggr)\|f\|_\infty
\nonumber\\[-8pt]\\[-8pt]
\eqntext{\mbox{for } x,y\in B_E(x_0,r/2).}
\end{eqnarray}

Now, for distinct $x,y\in D_n$ with $x\in B_E(0,c_1n^{1/2})$ and
$(\delta_n^{1/2})^2\le|x-y|< 1/4$ (note that since $|x-y|\ge\delta_n$
for distinct $x$ and $y$, the first inequality always hold), let
$x_0=x$ and $r=|x-y|^{1/2}$.
Then $\delta_n\le r< 1/2$ and $y\in B_E(x_0,r/2)$ (because
$|x_0-y|=r^2<r/2$). Thus, we can apply (\ref{cbe01}) to obtain
\begin{eqnarray*}
\bigl\llvert U_n^\lambda f(x)-U_n^\lambda
f(y)\bigr\rrvert&\leq& c \bigl(|x-y|^{\beta
/2}+\lambda^{-1}|x-y|^{\gamma/2}
\bigr)\|f\|_\infty
\\
&\leq& c\bigl(1+\lambda^{-1}\bigr)|x-y|^{(\beta\wedge\gamma)/2}\|f
\|_\infty.
\end{eqnarray*}
So, (\ref{cbe000}) holds with $C=c(1+\lambda^{-1})$ and $\gamma
'=(\beta
\wedge\gamma)/2$. The result at (\ref{cbe00110}) is immediate from
(\ref{cbe000}).
\end{pf}

We note that with an additional mild condition, we can further obtain
equi-H\"older continuity of the associated semigroup. (The next
proposition will only be used in the proof of Theorem~\ref{qip} below.)
Set $B_R:=B_E(0,R)\cap D_n$ for $R\in[2,\infty)$.
Denote by $X^{n, B_R}$ the subprocess of $X^n$ killed upon exiting $B_R$,
and $\{P_t^{n, B_R}; t\geq0\}$ the transition semigroup of $X^{n,
B_R}$. [When $R=\infty$, we set $(P_t^{n})_{t\geq0}:=(P_t^{n,
B_\infty
})_{t\geq0}$, i.e., the semigroup of $X^n$ itself.]
For $p\in[1,\infty]$, we use $\|\cdot\|_{p,n,R}$ to denote the
$L^p$-norm with respect to $m_n$ on $B_R$.

\begin{propn}\label{PR4} Let $R\in[2,\infty]$ and $t>0$. Suppose there exist
$c_1>0$ and $N_1\in\mathbb N$ (that may depend on $R$ and $t$) such
that for every $g\in L^1(B_R, m_n)$,
\[
{\bigl\llVert P^{n,B_R}_tg \bigr\rrVert}_{\infty,n,R}\leq
c_1{\llVert g \rrVert}_{1,n,R}\qquad\mbox{for all } n\ge
N_1. %
\]
Suppose in addition that Assumption~\ref{teoassumpHC} holds with
$X^{n,B_R}$ and $B_R$
in place of $X^n$ and $D_n$, respectively. It then holds that there
exist constants $c\in(0,\infty)$ and $N_2\geq1$ (that also may depend
on $R$ and $t$) such that
\[
\bigl\llvert P_t^{n, B_R} f(x)-P_t^{n, B_R}
f(y)\bigr\rrvert\leq c_2|x-y|^{\gamma
'}\|f\|_{2,n,R}
\]
for every $n\geq N_2$, $f\in L^2(B_R; m_n)$, and $m_n$-a.e. $x,y\in
B_{R/2}$ with $|x-y|< 1/4$. Here, $\gamma'$ is the constant of
Proposition~\ref{PR3}.
\end{propn}

\begin{pf} We follow \cite{BKK}, Proposition 3.4. For notational
simplicity, we
omit the superscripts $n$, $B_R$ on $P_t$
throughout the proof. Using spectral
representation theorem for self-adjoint operators,
there exist projection operators $E_\mu=E_\mu^{n,R}$ on the space
$L^2(B_R; m_n)$ such that
%
\begin{eqnarray}\label{C07}
f&=&\int_0^\infty dE_\mu(f),\qquad
P_tf=\int_0^\infty e^{-\mu t}
\,dE_\mu(f),
\nonumber\\[-8pt]\\[-8pt]
U^\lambda f&=&\int_0^\infty
\frac{1}{\lambda+\mu} \,dE_\mu(f).\nonumber
\end{eqnarray}
Define
\[
h=\int_0^\infty(\lambda+\mu)e^{-\mu t}
\,dE_\mu(f).
\]
Since $\sup_\mu(\lambda+\mu)^2 e^{-2\mu t}\leq c$, we have
\[
{\llVert h \rrVert}_2^2=\int_0^\infty(
\lambda+\mu)^2 e^{-2\mu t} \,d{ \bigl\langle E_\mu
(f), E_\mu(f) \bigr\rangle} \leq c\int_0^\infty
d{ \bigl\langle E_\mu(f), E_\mu(f) \bigr\rangle}=c{
\llVert f \rrVert}_2^2,
\]
where for $f,g\in L^2$, ${ \langle f,g \rangle}$ is the\vspace*{1pt}
inner product of $f$ and
$g$ in $L^2$. Thus, $h$ is a well-defined function in $L^2$.

Now,\vspace*{1pt} suppose $g\in L^1$. By the assumption, ${\llVert P_tg \rrVert
}_\infty\leq
c{\llVert g \rrVert}_1$, from which it follows that ${\llVert P_tg
\rrVert}_2\leq c{\llVert g \rrVert}_1$.
Since $\sup_\mu(\lambda+\mu) e^{-\mu t/2}\leq c$, using
Cauchy--Schwarz, we have
\begin{eqnarray*}
{ \langle h,g \rangle}&=&\int_0^\infty(\lambda+
\mu)e^{-\mu t} \,d{ \bigl\langle E_\mu(f), g \bigr\rangle}
\\
&\leq&\biggl(\int_0^\infty(\lambda+
\mu)e^{-\mu t} \,d \bigl\langle E_\mu(f), f\bigr\rangle\biggr)
^{1/2} \biggl(\int_0^\infty(\lambda+
\mu)e^{-\mu t} \,d{ \bigl\langle E_\mu(g), g \bigr\rangle}
\biggr)^{1/2}
\\
&\leq& c \biggl(\int_0^\infty d{ \bigl\langle
E_\mu(f), f \bigr\rangle} \biggr)^{1/2} \biggl(\int
_0^\infty e^{-\mu t/2} \,d{ \bigl\langle
E_\mu(g), g \bigr\rangle} \biggr)^{1/2}
\\
&=&c {\llVert f \rrVert}_2{\llVert P_{t/2}g \rrVert
}_2\leq c'{\llVert f \rrVert}_2 {\llVert g
\rrVert}_1.
\end{eqnarray*}
Taking the supremum over $g\in L^1$ with $L^1$ norm less than 1, this
yields ${\llVert h \rrVert}_\infty\leq c{\llVert f
\rrVert}_2$. Finally, by (\ref{C07}),
\[
U^\lambda h=\int_0^\infty
e^{-\mu t} \,dE_\mu(f)=P_tf\qquad\mbox{a.e.}
\]
and so the H\"older continuity of $P_tf$ follows from Proposition~\ref{PR3}.
\end{pf}

Let ${\mathbb D}(\mathbb{R}_+, \overline{D})$ be the space of right
continuous functions on $\mathbb{R}_+$ having left limits and
taking values in $\overline{D}$ that is equipped with
the Skorohod topology. For $t\geq0$, we use $X_t$ to denote
the coordinate projection map on ${\mathbb D}(\mathbb{R}_+, \overline{D})$;
that is, $X_t (\omega) = \omega(t)$ for $\omega\in{\mathbb
D}(\mathbb{R}_+,
\overline{D})$.
For subsequential convergence to a diffusion, we need the following.

\begin{assump}\label{teoassumpHC-2}
(i) For any sequence $x_n\rightarrow x$ with $x_n\in D_n$,
$\{\mathbb{P}^n_{x_n}\}$ is tight in ${\mathbb D}(\mathbb{R}_+,
\overline{D})$.

(ii)
For any sequence $x_n\rightarrow x$ with $x_n\in D_n$ and any
$\varepsilon>0$,
\[
\lim_{\delta\to0} \limsup_{n \to\infty}
\mathbb{P}^n_{x_n} \bigl( J\bigl(X^n,\delta
\bigr) >\varepsilon\bigr) =0,
\]
where
$J(X,\delta):=\int_0^\infty e^{-u} ( 1\wedge\sup_{\delta\le
t\le
u}|X_t-X_{t-\delta}| ) \,du$.
\end{assump}

We need the following well-known fact (see, e.g., \cite{baba89}, Lemma 6.4)
in the proof of Proposition~\ref{clttight}. For readers' convenience,
we provide a proof here.

\begin{lem}\label{BBlem64} Let
$K$ be a compact subset of $\mathbb R^d$.
Suppose $f$ and $f_k$, $k\in\mathbb N$, are functions on $K$ such that
$ \lim_{k\to\infty} f_k(y_k)=f(y)$ whenever $y_k\in K$ converges to~$y$. Then $f$ is continuous on $K$ and $f_k$ converges to $f$ uniformly
on $K$.
\end{lem}

\begin{pf} We first show that $f$ is continuous on $K$. Fix $x_0\in K$. Let
$x_k$ be any sequence in $K$ that converges to it.
Since $\lim_{i\to\infty} f_i (x) = f(x)$ for every $x\in K$,
there is a sequence $n_k\in{\mathbb N}$ that increases to infinity so that
$|f_{n_k}(x_k)-f(x_k)|\leq2^{-k}$ for every $k\geq1$. Since $\lim_{k\to
\infty}
f_{n_k}(x_k)=f(x_0)$, it follows that $\lim_{k\to\infty} |f(x_0)-f(x_k)|=0$.
This shows that $f$ is continuous at $x_0$, and hence on $K$.

We next show that $f_k$ converges uniformly to $f$ on $K$. Suppose not.
Then there is $\varepsilon>0$ so that for every $k\geq1$,
there are $n_k\geq k$ and $x_{n_k}\in K$ so that
$|f_{n_k}(x_{n_k})-f(x_{n_k})|>\varepsilon$.
Since $K$ is compact, by selecting a subsequence if necessary, we may
assume without loss of
generality that $x_{n_k} \to x_0\in K$.
As $\lim_{k\to\infty} f_{n_k}(x_{n_k}) =f(x_0)$ by the assumption,
we have $\liminf_{k\to\infty} |f(x_0)-f(x_{n_k})|\geq\varepsilon
$. This
contradicts to the
fact that $f$ is continuous on $K$.
\end{pf}

Now, applying the argument in \cite{baba89}, Section~6, we can prove that
any subsequential limit of the laws of $X^n$ under $\mathbb
{P}^n_{x_n}$ is the
law of a symmetric diffusion. For this, we need to introduce a
projection map from $\overline{D}$ to $D_n$. For each $n\geq1$, let
$\phi_n\dvtx
\overline{D}\to D_n$ be a map that projects each $x\in\overline{D}$ to
some $\phi
_n(x)\in D_n$ that minimizes $|x-y|$ over $y\in D_n$ (if there is more
than one such point that does this, we choose and fix one). If needed,
we extend a function $f$ defined on $D_n$ to be a function on
$\overline{D}$
by setting $f(x)= f(\phi_n(x))$. Note\vspace*{1pt} that each $D_n$ supports the
measure $m_n$ that converges vaguely to $m$. This implies that for each
$x\in\overline{D}$ and $r>0$, there is an $N\geq1$ so that $\phi_n
(x)\in
B_E(x, r)$ for every $n\geq N$.
From this, one concludes that
%
\begin{equation}
\label{e214} \qquad\phi_n (x_n) \to x_0\qquad\mbox{for every sequence } x_n\in\overline{D}\mbox{ that converges to } x_0.
\end{equation}

\begin{propn}\label{clttight}
Suppose that Assumptions~\ref{teoassumpHC} and~\ref{teoassumpHC-2}
hold and that $\{X^{n}, \mathbb{P}^n_x, x\in D_n\}$ is conservative for
sufficiently large $n$.
For every subsequence~$\{n_j\}$, there exists a sub-subsequence $\{
n_{j(k)}\}$ and
a continuous conservative $m$-symmetric Hunt process $(\widetilde X,
\widetilde{\mathbb{P}}
_x, x\in\overline{D})$ such that for every $x_{n_{j(k)}}\to x$,
$\mathbb{P}
^{n_{j(k)}}_{x_{n_j(k)}}$ converges weakly in $\mathbb{D}([0,\infty
),\overline{D})$ to $\widetilde{\mathbb{P}}^x$.
\end{propn}

\begin{pf} For notational simplicity, let us relabel the subsequence as
$\{n\}$.
We first claim that there exists\vspace*{-1pt} a (sub-)subsequence $\{n_j\}$ such
that $U_{n_j}^\lambda f$ converges uniformly on compact sets for each
$\lambda
>0$ and $f\in C_b(\overline{D})$.
Indeed, let $\{\lambda_i\}$ be a dense subset of $(0,\infty)$
and $\{f_k\}$ a sequence of
functions in $C_b(\overline{D})$
such that $\|f_k\|_\infty\le1$ and whose linear span is dense in
$(C_b (\overline{D}), \| \cdot\|_\infty)$.
For fixed $m$ and $i$,
by Proposition~\ref{PR3} and the Ascoli--Arzel\`a theorem, there is a
subsequence of $U_n^{\lambda_i}f_k$ that
converges uniformly on compact sets.
By a diagonal selection procedure,
we can choose a subsequence $\{n_j\}$ such that $U_{n_j}^{\lambda_i}f_k$
converges uniformly on compact sets for every $m$ and $i$
to a H\"older continuous function which we denote as $U^{\lambda_i}f_k$.
Noting that\vspace*{-2pt}
%
\begin{eqnarray}\label{eqresol25}
U_n^\lambda-U_n^\beta&=&(\beta-\lambda)U_n^\lambda U_n^\beta,\nonumber
\\[-2pt]
\bigl\|U_n^\lambda\bigr\|_{\infty\to\infty}&\le&\frac{1}{\lambda},
\\[-2pt]
\bigl\| U_n^\lambda-U_n^\beta
\bigr\|_{\infty\to\infty}&\le&\frac{\beta-\lambda}{\lambda
\beta},\nonumber
\end{eqnarray}
a\vspace*{-1pt} careful limiting argument shows
that $U_{n_j}^{\lambda}f$ converges uniformly on compact sets, say to
$U^{\lambda}f$, for any $\lambda>0$ and any continuous function
$f$, and
(\ref{eqresol25}) holds as well for $\{U^\lambda\}$.
By the equi-continuity of $U_{n_j}^{\lambda}f$, we also have
$U_{n_j}^{\lambda}f(x_{n_j})\to U^\lambda f(x)$
for each $x_{n_j}\in D_{n_j}$ that converges to $x\in\overline{D}$.

We\vspace*{-2pt} next claim that $\mathbb{P}^{n_j}_{x_{n_j}}$ converges weakly, say
to $\widetilde{\mathbb{P}}_x$. Indeed, by Assumption~\ref{teoassumpHC-2}(i),
$\{\mathbb{P}^{n_j}_{x_{n_j}}\}$ is tight,
so it suffices to show that any two limit points agree.
Let $\mathbb{P}'$ and $\mathbb{P}''$ be any two limit points. Then
one sees that
\[
\mathbb{E}' \biggl[ \int_0^\infty
e^{-\lambda s}f(X_s)\,ds \biggr]=U^\lambda f(x)=\mathbb{E}
'' \biggl[\int_0^\infty
e^{-\lambda s}f(X_s)\,ds \biggr]
\]
for any $f\in C_b(\overline{D})$.
So, by the uniqueness of the Laplace transform,
\[
\mathbb{E}'\bigl[f(X_s)\bigr]= \mathbb{E}''
\bigl[ f(X_s) \bigr]
\]
for almost all $s\ge0$
and hence for every $s\geq0$ since $s\to X_s$ is right continuous.
So, the one-dimensional
distributions of $X_t$ under $\mathbb{P}'$ and $\mathbb{P}''$ are
the same. Set
$P_sf(x):=\mathbb{E}'f(X_s)$. We have $P_s^{n_j}f(x_{n_j})\to P_sf(x)$
for every sequence $x_{n_j}\in D_{n_j}$ that converges to $x$.
Recall the restriction map $\phi_n$ introduced proceeding the statement
of this theorem. It follows from (\ref{e214}) that
$P_s^{n_j}(f\circ\phi_{n_j}) (y_{n_j})\to P_sf(y)$ for every sequence
$y_{n_j}\in\overline{D}$ that converges to $y$.
Thus, by Lemma~\ref{BBlem64},
$P^{n_j}_s (f\circ\phi_{n_j})$ converges to $P_s f$ uniformly on
compact subsets of $\overline{D}$ and $P_s f \in C_b (\overline{D})$ for every
$f\in C_c(\overline{D})$.
For $f, g\in C_c (\overline{D})$
and $0\le s<t$, by the\vspace*{1pt} Markov property of
$X$ under~$\mathbb{P}^{n_j}_{x_{n_j}}$,
\begin{eqnarray*}
\mathbb{E}^{n_j}_{x_{n_j}} \bigl[g(X_s)f(X_t)
\bigr] &=&\mathbb{E}^{n_j}_{x_{n_j}} \bigl[\bigl(
\bigl(P_{t-s}^{n}f\bigr)g\bigr) (X_s) \bigr]
\\[-2pt]
&=&\mathbb{E}^{n_j}_{x_{n_j}} \bigl[\bigl((P_{t-s}f)g
\bigr) (X_s) \bigr] +\mathbb{E}^{n_j}_{x_{n_j}}
\bigl[\bigl(\bigl(P_{t-s}^{n}f-P_{t-s}f\bigr)g
\bigr) (X_s) \bigr].
\end{eqnarray*}
The first term of the right-hand side converges to $\mathbb{E}
'[((P_{t-s}f)g)(X_s)]$ by the above proof, while the second term goes
to $0$ since $P_{t-s}^{n_j}f\to P_{t-s}f$ uniformly on compact sets.
Repeating this, we conclude that for every $k\geq1$ and every $0< s_1<s_2
<\cdots< s_k$ and $f_j \in C_c(\overline{D})$,
%
\begin{eqnarray}\label{e210}
\mathbb{E}' \Biggl[ \prod
_{j=1}^k f_j (X_{s_j})
\Biggr] &=& \mathbb{E}'' \Biggl[ \prod
_{j=1}^k f_j (X_{s_j})
\Biggr]
\nonumber\\[-8pt]\\[-8pt]
&=& P_{s_1} \bigl( f_1 P_{s_2-s_1}
\bigl(f_2 P_{s_3-s_2} (f_3 \cdots)\bigr) \bigr)
(x).\nonumber
\end{eqnarray}
This proves that the finite-dimensional distributions of $X$ under
$\mathbb{P}'$ and $\mathbb{P}''$ are the same. Consequently, $\mathbb
{P}'=\mathbb{P}''$, which we
now denote as $\widetilde{\mathbb{P}}_x$.
Moreover, (\ref{e210}) shows that $(X, \widetilde{\mathbb{P}}_x,
x\in\overline{D})$
is a Markov process with transition semigroup $\{P_t, t\geq0\}$.

Next, we show that $\{\widetilde{\mathbb{P}}_x\dvtx x\in\overline{D}\}$
is a strong Markov process.
Note that
$X$ is conservative with $\widetilde{\mathbb{P}}_x(X_0=x)=1$, and
under Assumption~\ref{teoassumpHC-2}(ii), $X_t$ is continuous a.s.
under $\widetilde{\mathbb{P}}_x$.
We also have $P_t f \in C_b(\overline{D})$ for $f\in C_c (\overline{D})$.
It is easy to deduce from these properties and (\ref{e210})
that for every $f\in C_c (\overline{D})$ and every stopping time $T$,
\[
\mathbb{E}_x \bigl[ f(X_{T+t})| \mathcal{F}_{T+}
\bigr] = P_t f (X_T), \qquad x\in\overline{D}.
\]
See the proof of Theorem 2.3.1 on page~56 of \cite{CW}.
From it, one gets the strong Markov property of $X$ by a standard
measure-theoretic argument (see page~57 of~\cite{CW}). Since $X$ is
continuous and has infinite lifetime, this in fact shows that $X$ is a
continuous conservative Hunt process.

Finally, for $f, g\in C_c (\overline{D})$, by the convergence of
semigroups and vague convergence of measures,
it holds that, for every $t>0$,
\[
\int_{D_n} \bigl(P_t^{n_j}f\bigr) (x)
g(x)m_{n_j}(dx)\rightarrow\int_{\overline{D}}
(P_tf) (x) g(x)m(dx).
\]
Since $X^{n_j}$ is $m_n$-symmetric, this readily yields the desired
$m$-symmetry of $\widetilde{X}$.
\end{pf}

\begin{rem}
Note that we did not use any special properties of the
Euclidean metric in this section, so that all the arguments in this
section can be extended to a metric measure space without any changes.
\end{rem}

By Theorem~\ref{teosilvs}, Remark~\ref{teoconserrem}(i) and
Proposition~\ref{clttight}, we see that in order to prove
$(X^{n}, \mathbb{P}^n_{x_n})$ converges weakly to $(X, \mathbb{P}_x)$ in
$\mathbb{D}([0,\infty),\overline{D})$ as $n\to\infty$, it suffices
to verify
conditions (ii) and (iii) in Theorem~\ref{teosilvs},
vague convergence of the measure $m_n$ to $m$ on $\overline{D}$,
Assumptions~\ref{teoassumpHC} and~\ref{teoassumpHC-2}, and the
conservativeness of $(X^{n}, \mathbb{P}^n_{x_n})$ for each $n\in
\mathbb{N}$.

\section{Random conductance model in unbounded domains}\label{sec3}\label{percsec}

In this section, we will obtain, as a first application of our theorem,
a quenched invariance principle for random walk among random
conductances on half-spaces, quarter-spaces, etc. The assumptions we
make on the random conductances include the supercritical percolation
model, and random conductances boun\-ded uniformly from below and with
finite first moments. For the main conclusion, see Theorem~\ref{RCMresult}.

Fix $d_1,d_2\in\mathbb{Z}_+$ such that $d_1\geq1$ and
$d:=d_1+d_2\geq2$.
Define a graph $(\mathbb{L},E_{\mathbb{L}})$ by setting $\mathbb
{L}:=\mathbb{Z}^{d_1}_+\times\mathbb{Z}^{d_2}$ and $E_{\mathbb
{L}}:=\{
e=\{x,y\}\dvtx x,y\in\mathbb{L}, |x-y|=1\}$. Given $\mathcal
{O}\subseteq
E_\mathbb{L}$, let $\mathcal{C}_\infty(\mathbb{L},\mathcal{O})$ be the
infinite connected cluster of $(\mathbb{L},\mathcal{O})$, provided it
exists and is unique [otherwise set $\mathcal{C}_\infty(\mathbb
{L},\mathcal{O}):=\varnothing$].

Let
$\mu=(\mu_e)_{e\in E_\mathbb{L}}$ be a collection of independent and
identically distributed random variables on $[0,\infty)$, defined on a
probability space
$(\Omega, \mathbb{P})$ such that
%
\begin{equation}
\label{mucond1} p_1:=\mathbb{P} (\mu_e>0
)>p_c^{\mathrm{bond}}\bigl(\mathbb{Z}^d\bigr),
\end{equation}
where $p_c^{\mathrm{bond}}(\mathbb{Z}^d)\in(0,1)$ is the critical
probability for bond percolation on $\mathbb{Z}^d$.
We assume that there is $c>0$ so that
%
\begin{equation}
\label{mucond2} \mathbb{P} \bigl(\mu_e\in(0,c) \bigr)=0
\end{equation}
and
%
\begin{equation}
\label{mucond3} \mathbb{E} (\mu_e )<\infty.
\end{equation}
This
framework includes the special cases of supercritical percolation
[where
$\mathbb{P} (\mu_e=1 )=p_1=1-\mathbb{P} (\mu_e=0
)$]
and the random conductance model
with conductances bounded from below [i.e., $\mathbb{P} (c\le
\mu
_e<\infty)=1$ for some $c>0$] and having finite first moments.
For each $x\in\mathbb{L}$, set $\mu_x=\sum_{y\sim x}\mu_{xy}$. Set
\[
\mathcal{O}_1:= \{e\in E_{\mathbb{L}}\dvtx\mu_e>0 \},
\qquad\mathcal{C}_1:=\mathcal{C}_\infty(\mathbb{L},
\mathcal{O}_1).
\]
Note that for the random conductance model bounded from below,
$\mathcal
{C}_1=\mathbb{L}$.

For each realization of $\mathcal{C}_1$,
there is
a continuous time Markov chain $Y=(Y_t)_{t\geq0}$ on $\mathcal{C}_1$
with transition probabilities $P(x,y)=\mu_{xy}/\mu_x$, and the holding
time at
each $x\in\mathcal{C}_1$
being the exponential distribution with mean $\mu_x^{-1}$.
Such a Markov chain is sometimes called a variable speed random walk
(VSRW). The corresponding Dirichlet form is $(\mathcal{E},
L^2(\mathcal{C}_1;
\nu))$, where $\nu$ is
the counting
measure on $\mathcal{C}_1$
and
\[
\mathcal{E}(f,g)=\frac{1}2\sum_{x,y\in\mathcal{C}_1, {x\sim
y}}
\bigl(f(x)-f(y)\bigr) \bigl(g(x)-g(y)\bigr)\mu_{xy} \qquad\mbox{for }
f,g\in L^2(\mathcal{C}_1; \nu).
\]
The corresponding discrete Laplace operator is
${\mathcal L}_V f(x)=\sum_y(f(y)-f(x))\mu_{xy}$. For each $f,g$ that
have finite support, we have
\[
{\mathcal E}(f,g)=-\sum_{x\in\mathcal{C}_1}({\mathcal
L}_Vf) (x)g(x).
\]
We will establish a quenched invariance principle for $Y$ in
Section~\ref{qipsec}, but we first need to derive some preliminary
estimates regarding the geometry of $\mathcal{C}_1$ and the heat kernel
associated with $Y$.

\subsection{Percolation estimates}\label{sec3.1}\label{percosect}

In this section, we derive a number of useful properties of the
underlying percolation cluster $\mathcal{C}_1$. Most importantly, we
introduce the concept of ``good'' and ``very good'' balls for the model
and provide estimates for the probability of such occurring; see
Definition~\ref{good} and Proposition~\ref{teoVNG-lem} below.

Since a variety of percolation models will appear in the course of
this paper, let us now make explicit that
the critical probability for bond/site percolation on an infinite
connected graph containing the vertex $0$ is
\begin{eqnarray*}
p_c &:=& \inf\bigl\{p\in[0,1]\dvtx \mathbb{P}_p (\mbox{$0$ is
in an infinite connected}
\\
&&\hspace*{87pt} \mbox{cluster of open bonds/sites} )>0 \bigr\},
\end{eqnarray*}
where $\mathbb{P}_p$ is the law of parameter $p$ bond/site percolation
on the graph in question. Note in particular that the critical
probability for bond percolation on $\mathbb{L}$ is identical to
$p_c^{\mathrm{bond}}(\mathbb{Z}^d)$; see \cite{grim}, Theorem 7.2
(which is the bond percolation version of a result originally proved as
\cite{GM}, Theorem~A). Recall
\[
\mathcal{O}_1:= \{e\in E_{\mathbb{L}}\dvtx\mu_e>0 \},
\qquad\mathcal{C}_1:=\mathcal{C}_\infty(\mathbb{L},
\mathcal{O}_1).
\]
That $\mathcal{C}_1$ is nonempty almost-surely is guaranteed by \cite
{BGN}, Corollary to Theorem 1.1 (this covers $d\geq3$, and, as is
commented there, the case $d=2$ can be tackled using techniques from
\cite{Harris}).

Now, suppose that $\mu$ is actually a restriction of independent and
identically distributed (under $\mathbb{P})$ random variables $(\mu
_e)_{e\in E_{\mathbb{Z}^d}}$, where $E_{\mathbb{Z}^d}$ are the usual
nearest-neighbor edges for the integer lattice $\mathbb{Z}^d$, and define
\[
\widetilde{\mathcal{O}}_1:= \{e\in E_{\mathbb{Z}^d}\dvtx
\mu_e>0 \},\qquad\widetilde{\mathcal{C}}_1:=
\mathcal{C}_\infty\bigl(\mathbb{Z}^d,\widetilde{\mathcal{O}}_1\bigr).
\]
For sufficiently large $K$, so that
\[
q=q(K):=\mathbb{P}\bigl(0<\mu_e<K^{-1}\bigr)+\mathbb{P}(
\mu_e>K)<p_1-p_c^{\mathrm
{bond}}\bigl(
\mathbb{Z}^d\bigr)
\]
and writing $\widetilde{\mathcal{O}}_I:=\{e\in E_{\mathbb{Z}^d}\dvtx\mu
_e\in I\}
$ for $I\subseteq[0,\infty)$, we let
\begin{eqnarray}
\widetilde{\mathcal{O}}_R&:=&\widetilde{\mathcal{O}}_{(0,K^{-1})\cup
(K,\infty
)},
\nonumber
\\
\widetilde{\mathcal{O}}_S&:=& \bigl\{e\in\widetilde{\mathcal{O}}_1\dvtx e
\cap e'\neq\varnothing\mbox{ for some }e'\in\widetilde{\mathcal{O}}_R
\bigr\},
\nonumber
\\
\widetilde{\mathcal{O}}_2&:=&\widetilde{\mathcal{O}}_1\setminus
\widetilde{\mathcal{O}}_S.
\nonumber
\end{eqnarray}
We will also define $\mathcal{O}_2:=\widetilde{\mathcal{O}}_2\cap
E_{\mathbb
{L}}$, and set $\mathcal{C}_2:=\mathcal{C}_\infty(\mathbb
{L},\mathcal
{O}_2)$---the next lemma will guarantee that this set is nonempty
almost-surely.

To represent the set of ``holes,'' let $\mathcal{H}:=\mathcal
{C}_1\setminus\mathcal{C}_2$. Moreover, for $x\in\mathcal{C}_1$, let
$\mathcal{H}(x)$ be the connected component of $\mathcal
{C}_1\setminus
\mathcal{C}_2$ containing $x$. The following lemma provides control on
the size of these components. Since its proof is a somewhat technical
adaptation to our setting of that
used to establish \cite{ABDH}, Lemma 2.3, which dealt with the whole
$\mathbb{Z}^d$ model, we defer this to the \hyperref[sec5]{Appendix}. Note, though, that
in the percolation case (i.e., when $\mu_e$ are Bernoulli random
variables) or the uniformly elliptic random conductor case, the proof
of the result is immediate; indeed, for large enough $K$, we have that
$\mathcal{O}_2=\mathcal{O}_1$, and so $\mathcal{H}=\varnothing$.

\begin{lem}\label{holelem}
For sufficiently large $K$,
the following hold:
\begin{longlist}[(ii)]
\item[(i)] All the connected components of $\mathcal{H}$ are finite.
Furthermore, there exist constants $c_1,c_2$ such that: for each $x\in
\mathbb{L}$,
\[
\mathbb{P} \bigl(x\in\mathcal{C}_1\mbox{ and } \operatorname{diam}\bigl(
\mathcal{H}(x)\bigr)\geq n \bigr)\leq c_1e^{-c_2n},
\]
where $\operatorname{diam}$ denotes the diameter with respect to the
$\ell
_\infty
$ metric on $\mathbb{Z}^d$.

\item[(ii)] There exists a constant $\alpha$ such that, $\mathbb{P}$-a.s., for
large enough $n$, the volume of any hole intersecting the box
$[-n,n]^d\cap\mathbb{L}$ is bounded above by $(\log n)^\alpha$.
\end{longlist}
\end{lem}

In what follows, we will need to make comparisons between two graph
metrics on $(\mathcal{C}_1,\mathcal{O}_1)$, and the Euclidean metric.
The first of these, $d_1$, will simply be defined to be the shortest
path metric on $(\mathcal{C}_1,\mathcal{O}_1)$, considered as an
unweighted graph. To define the second metric, $\bar{d}_1$, we follow
\cite{BarDeu08} by defining edge weights
\[
t(e):=C_A\wedge\mu_e^{-1/2},
\]
where $C_A<\infty$ is a deterministic constant, and then letting
$\bar
{d}_1$ be the shortest path metric on $(\mathcal{C}_1,\mathcal{O}_1)$,
considered\vspace*{1pt} as a weighted graph (in \cite{BarDeu08}, the analogous
metric was denoted $\tilde{d}$). We note that the latter metric on
$\mathcal{C}_1$ satisfies
\[
\bigl(C_A^{-2}\vee\mu_{\{y,z\}} \bigr)\bigl\llvert
\bar{d}_1(x,y)-\bar{d}_1(x,z)\bigr\rrvert
^2\leq1
\]
for any $x,y,z\in\mathcal{C}_1$ with $\{y,z\}\in\mathcal{O}_1$. Observe
that, since the weights $t(e)$ are bounded above by $C_A$, we
immediately have that $\bar{d}_1$ is bounded above by $C_A d_1$, Hence,
the following lemma establishes both $\bar{d}_1$ and $d_1$ are
comparable to the Euclidean one. An easy consequence of this is the
comparability of balls in the different metrics; see Lemma~\ref{ballcomp}. The proofs of both these results are deferred to the \hyperref[sec5]{Appendix}.

\begin{lem}\label{distlem} There exist constants $c_1,c_2,c_3$ such
that: for
$R\geq1$,
%
\begin{equation}
\label{upperd1} \mathop{\sup_{x,y\in\mathbb{L}\dvtx}}_{|x-y|\leq
R}\mathbb{P}
\bigl(x,y\in\mathcal{C}_1\mbox{ and } {d}_1(x,y)\geq
c_1R \bigr)\leq c_2e^{-c_3R}
\end{equation}
and also, for every $x,y\in\mathbb{L}$,
%
\begin{equation}
\label{lowerd1} \mathbb{P} \bigl(x,y\in\mathcal{C}_1\mbox{ and }
\bar{d}_1(x,y)\leq c_1^{-1}|x-y| \bigr)
\leq c_2e^{-c_3|x-y|}.
\end{equation}
\end{lem}

\begin{lem}\label{ballcomp} There exist constants $c_1,c_2,c_3,c_4$ such
that: for every $x\in\mathbb{L}$, $R\geq1$,
\begin{eqnarray*}
&& \mathbb{P} \bigl( \{x\in\mathcal{C}_1 \}\cap\bigl\{ \mathcal
{C}_1\cap B_E(x,c_1R)\subseteq
B_1(x,R)\subseteq\overline{B}_1(x,C_AR)
\subseteq B_E(x,c_2R) \bigr\}^c \bigr)
\\
&&\qquad \leq c_3e^{-c_4R},
\end{eqnarray*}
where $B_1(x,R)$ is a ball in the metric space $(\mathcal{C}_1,d_1)$,
$\overline{B}_1(x,R)$ is a ball in $(\mathcal{C}_1,\bar
{d}_1)$, and
$B_E(x,R)$ is a Euclidean ball.
\end{lem}

We continue by adapting a definition for ``good'' and ``very good''
balls from \cite{ABDH}. In preparation for this, we define ${\mu
}^0_e:=\mathbf{1}_{\{e\in{\mathcal{O}}_1\}}$, set ${\mu}^0_x:=\sum_{y\in
\mathbb{L}}{\mu}_{\{x,y\}}^0$ for $x\in\mathbb{L}$, and then
extend ${\mu}^0$ to a measure on $\mathbb{L}$. Moreover, we set
$\beta:=1-2(1+d)^{-1}$.

\begin{defn}\label{good} (i) Let $C_V, C_P, C_W, C_R, C_D$ be fixed strictly
positive constants. We say the pair $(x,R)\in\mathcal{C}_1\times
\mathbb
{R}_+$ is \emph{good} if
%
\begin{eqnarray}
\label{gp1} {B}_1\bigl(x,C_A^{-1}r\bigr)
&\subseteq& \overline{B}_1(x,r)\subseteq B_1(x,C_Dr)\qquad\forall r\geq R,
\\
\label{gp2} |y-z|&\geq& C_R^{-1}R\qquad\forall y\in
\overline{B}_1(x,R/2), z\in\overline{B}_1(x,8R/9)^c,
\\
\label{gp3} C_VR^d&\leq&\mu^0\bigl(
\overline{B}_1(x,R)\bigr),
\\
\label{gp4} \operatorname{diam} \bigl(\mathcal{H}(y)\bigr)&\le& R^\beta
\qquad\forall y\in B_E(x,R)\cap\mathcal{C}_1
\end{eqnarray}
and the weak Poincar\'{e} inequality
%
\begin{equation}
\label{poincare}\qquad \sum_{y\in{B}_1(x,R)} \bigl(f(y)-
\check{f}_{{B}_1(x,R)} \bigr)^2\mu_y^0
\leq C_P R^2\mathop{\sum_{y,z\in{B}_1(x,C_WR)\dvtx}}_{\{y,z\}\in
\mathcal
{O}_1}
\bigl\llvert f(y)-f(z)\bigr\rrvert^2
\end{equation}
holds for every $f\dvtx{B}_1(x,C_WR)\rightarrow\mathbb{R}$. [Here,
$\check
{f}_{{B}_1(x,R)}$ is the value which minimizes the left-hand side of
(\ref{poincare}).]

(ii) We say a pair $(x,R)\in\mathcal{C}_1\times\mathbb{R}_+$ is
\emph{very good} if: there exists $N=N_{(x,R)}$ such that $(y,r)$ is good
whenever $y\in\overline{B}_1(x,R)$ and $N\leq r\leq R$. We can always
assume that $N\geq2$. Moreover, if $N\leq M$, we will say that $(x,R)$
is $M$-very good.

(iii) Let $\alpha\in(0,1]$. For $x\in\mathcal{C}_1$, we define
$R_x^{(\alpha)}$ to be the smallest integer $M$ such that $(x,R)$ is
$R^\alpha$-very good for all $R\geq M$. We set $R_x^{(\alpha)}=0$ if
$x\notin\mathcal{C}_1$.
\end{defn}

The following proposition, which is an adaptation of \cite{ABDH},
Proposition 2.8, provides bounds for the probabilities of these events and
for the distribution of $R_x^{(\alpha)}$.

\begin{propn}\label{teoVNG-lem} There exist $c_1,c_2,C_V, C_P, C_W, C_R,
C_D$ (depending on the law of $\mu$ and the dimension $d$) such that
the following holds. For $x \in\mathbb{L}$, $R\geq1$, $\alpha\in
(0, 1]$,
%
\begin{eqnarray}
\label{a1} \mathbb{P} \bigl(x\in\mathcal{C}_1, (x,R)\mbox{ is not
good} \bigr)&\leq& c_1e^{-c_2R^\beta},
\\
\label{a2} \mathbb{P} \bigl(x\in\mathcal{C}_1, (x,R)\mbox{ is not
$R^\alpha$-very good} \bigr)&\leq& c_1e^{-c_2R^{\alpha\beta}}.
\end{eqnarray}
Hence,
%
\begin{equation}
\label{a3} \mathbb{P} \bigl(x\in\mathcal{C}_1,
R_x^{(\alpha)}\geq n \bigr)\leq c_1e^{-c_2n^{\alpha\beta}}.
\end{equation}
\end{propn}

\begin{pf} That
\[
\mathbb{P} \bigl(x\in\mathcal{C}_1, \mbox{ (\ref{gp1}) does not
hold} \bigr)\leq c_3e^{-c_4R}
\]
is a straightforward consequence of Lemma~\ref{ballcomp}.

For the second property, we have
\begin{eqnarray*}
&& \mathbb{P} \bigl(x\in\mathcal{C}_1, \mbox{ (\ref{gp2})
does not hold} \bigr)
\\
&&\qquad =\mathbb{P} \bigl(x\in\mathcal{C}_1, \exists y\in\overline
{B}_1(x,R/2), z\in\overline{B}_1(x,8R/9)^c\dvtx
|y-z|< C_R^{-1}R \bigr)
\\
&&\qquad \leq\mathbb{P} \bigl(x\in\mathcal{C}_1, \exists y
\in{B}_E(x,c_5R)\cap\overline{B}_1(x,R/2),
\\
&&\hspace*{63pt} z\in\overline{B}_1(x,8R/9)^c\dvtx |y-z|<C_R^{-1}R
\bigr)
+c_6e^{-c_7R}
\\
&&\qquad \leq c_6e^{-c_7R}+\sum_{y\in{B}_E(x,c_5R)}
\sum_{z\in
B_E(y,C_R^{-1}R)}\mathbb{P} \bigl(y,z\in
\mathcal{C}_1, \bar{d}_1(y,z)>R/3 \bigr)
\\
&&\qquad \leq c_8e^{-c_9R},
\end{eqnarray*}
where we apply Lemma~\ref{ballcomp} to deduce the first inequality, and
(\ref{upperd1}) to obtain the final one.

For (\ref{gp3}), applying Lemma~\ref{ballcomp} again yields
\begin{eqnarray*}
&& \mathbb{P} \bigl(x\in\mathcal{C}_1, \mbox{ (\ref{gp3})
does not hold} \bigr)
\\
&&\qquad =\mathbb{P} \bigl(x\in\mathcal{C}_1, \mu^0\bigl(
\overline{B}_1(x,R)\bigr)<C_VR^d \bigr)
\\
&&\qquad \leq c_{10}e^{-c_{11}R}+\mathbb{P} \bigl(x\in
\mathcal{C}_1, |\mathcal{C}_1\cap{B}_E(x,c_{12}R)|<C_VR^d
\bigr).
\end{eqnarray*}
Now, let $Q\subseteq B_E(x,c_{12}R)\cap\mathbb{L}$ be a cube of
side-length $c_{13}R$ such that $\inf_{y\in\mathbb{L}\setminus
Q}|x-y|\geq c_{13} R/2$. Moreover, if we let $\mathcal{C}^+(Q)$ be the
largest connected component of the graph $(Q,\mathcal{O}_1)$, then
\begin{eqnarray*}
&& \mathbb{P} \bigl(x\in\mathcal{C}_1, \mbox{ (\ref{gp3})
does not hold} \bigr)
\\
&&\qquad \leq \mathbb{P} \bigl(x\in\mathcal{C}_1\cap\mathcal{C}^+(Q), |
\mathcal{C}_1\cap Q|<C_VR^d \bigr)+
\mathbb{P} \bigl(x\in\mathcal{C}_1\setminus\mathcal{C}^+(Q)
\bigr)+c_{10}e^{-c_{11}R}
\\
&&\qquad \leq\mathbb{P} \bigl(\bigl|\mathcal{C}^+(Q)\bigr|<C_VR^d
\bigr)+\mathbb{P} \bigl(x\in\widetilde{\mathcal{C}}_1\setminus
\mathcal{C}^+(Q) \bigr)+c_{10}e^{-c_{11}R}.
\end{eqnarray*}
This bound is now expressed in terms of the full $\mathbb{Z}^d$ model,
for which appropriate estimates already exist. In particular, the first
term here is bounded above by $\mathbb{P}(G(Q)^c)$, where $G(Q)$ is the
event that $|\mathcal{C}^+(Q)|\geq\frac{1}{2}\mathbb{P}(0\in\widetilde
{\mathcal{C}}_1)|Q|$
[recall that $\widetilde{\mathcal{C}}_1:=\mathcal{C}_\infty(\mathbb
{Z}^d,\widetilde{\mathcal{O}}_1)$]. Consequently, by simply translating the
relevant part of \cite{ABDH}, Lemma 2.6, to our setting (taking
$K=\infty
$), we obtain that it is bounded above by $c_{13}e^{-c_{14}R^\beta}$.
That the second term is bounded above by $c_{15}e^{-c_{16}R}$ can be
established by applying \cite{bar00}, Lemma 2.8.

To check the fourth property, we simply note
\[
{\mathbb{P} \bigl(x\in\mathcal{C}_1, \mbox{ (\ref{gp4}) does
not hold} \bigr)}\leq\sum_{y\in B_E(x,R)\cap\mathbb{L}}\mathbb{P} \bigl(y
\in\mathcal{C}_1, \operatorname{diam} \bigl(\mathcal{H}(y)\bigr)>
R^\beta\bigr),
\]
which may be bounded above by $c_{17}e^{-c_{18}R^\beta}$ by applying
Lemma~\ref{holelem}.

Finally, for the Poincar\'{e} inequality, we will apply \cite{bar00},
Proposition 2.12. In particular, this result yields that if
$Q$ is a cube of side-length $2R$ contained in $\mathbb{L}$, $\mathcal
{C}^+(Q)$ is the largest connected component of the graph $(Q,\mathcal
{O}_1)$, and $H(Q)$ is the event that
\[
\min_{a}\sum_{y\in\mathcal{C}^+(Q)}
\bigl(f(y)-a \bigr)^2\mu_y^0 \leq
CR^2\mathop{\sum_{y,z\in\mathcal{C}^+(Q)\dvtx}}_{\{y,z\}\in
\mathcal
{O}_1}
\bigl\llvert f(y)-f(z)\bigr\rrvert^2
\]
for every $f\dvtx\mathcal{C}^+(Q)\rightarrow\mathbb{R}$, then $\mathbb
{P}(H(Q)^c)\leq c_{19}e^{-c_{20}R^\beta}$. Furthermore, it is clear
that if $H(Q)$ holds and also $B_1(x,R)\subseteq\mathcal
{C}^+(Q)\subseteq B_1(x,c_{21}R)$, then (\ref{poincare}) holds. This
means that
\begin{eqnarray*}
&& \mathbb{P} \bigl(x\in\mathcal{C}_1, \mbox{ (\ref{poincare})
does not hold} \bigr)
\\
&&\qquad \leq c_{22}e^{-c_{23}R^\beta} +\mathbb{P} \bigl( \{x\in\mathcal
{C}_1 \}\cap\bigl\{B_1(x,R)\subseteq\mathcal{C}^+(Q)
\subseteq B_1(x,c_{21}R) \bigr\}^c \bigr),
\end{eqnarray*}
where $Q$ is chosen such that $B_E(x,R)\cap\mathbb{L}\subseteq Q$.
Noting as above that $\mathbb{P}(x\in\mathcal{C}_1\setminus
\mathcal
{C}^+(Q))\leq c_{24}e^{-c_{25}R}$, at the expense of adjusting
constants, we may replace $\{x\in\mathcal{C}_1\}$ by $\{x\in\mathcal
{C}_1\cap\mathcal{C}^+(Q)\}$ in the above bound. On the event $\{x\in
\mathcal{C}_1\cap\mathcal{C}^+(Q)\}\cap\{B_1(x,R)\subseteq\mathcal
{C}^+(Q)\}^c$, it is elementary to check that $B_1(x,R)\nsubseteq
B_E(x,R)$, which is impossible. Since $\mathcal{C}^+(Q)\subseteq
B_E(x,2R)$, we have thus shown that
\begin{eqnarray*}
&& \mathbb{P} \bigl(x\in\mathcal{C}_1, \mbox{ (\ref{poincare})
does not hold} \bigr)
\\
&&\qquad \leq c_{26}e^{-c_{27}R^\beta} +\mathbb{P} \bigl( \{x\in\mathcal
{C}_1 \}\cap\bigl\{B_E(x,2R)\subseteq
B_1(x,cR) \bigr\} ^c \bigr).
\end{eqnarray*}
By applying Lemma~\ref{ballcomp} once again, this expression is bounded
above by $c_{28}e^{-c_{29}R^\beta}$, and so we have completed the proof
of (\ref{a1}).

Given (\ref{a1}), a simple union bound subsequently yields (\ref{a2}),
and the inequality at (\ref{a3}) is a straightforward consequence of this.
\end{pf}

\begin{rem}\label{teoremddtil}
It only requires a simple argument to check that if $(x,R)$ is good and
$y\in\mathcal{C}_1$ satisfies $d_1(x,y)\geq C_DR$, then
\[
C_D^{-1}d_1(x,y)\leq\bar
{d}_1(x,y)\leq C_Ad_1(x,y)
\]
(cf. \cite{BarDeu08}, Lemma 2.10(a)).
\end{rem}

Finally, we state a bound that allows us to compare $\nu$ with $\tilde
{\nu}$, which is the measure-defined similarly from the whole $\mathbb
{Z}^d$ model, that is, uniform measure on~$\widetilde{\mathcal{C}}_1$. Its
proof can be found in the \hyperref[sec5]{Appendix}.

\begin{lem}\label{meascomp} There exists a constant $c$ such that if
$Q\subseteq\mathbb{L}$ is a cube of side length $n$, then
\[
\mathbb{P} \bigl(\tilde{\nu}(Q)-{\nu}(Q)\geq n^{d-1}(\log
n)^{d+1} \bigr)\leq cn^{-2}.
\]
\end{lem}

\subsection{Heat kernel estimates}\label{sec3.2}\label{hkestsec}

Let $D_n=n^{-1}\mathcal{C}_1$, $\overline{D}=\mathbb{R}^{d_1}_+\times
\mathbb{R}^{d_2}$
(recall $\mathbb{R}_+:=[0,\infty)$).
Let $Y$ be the VSRW on $\mathcal{C}_1$, and for a given realization of
$\mathcal{C}_1=\mathcal{C}_1(\omega)$, $\omega\in\Omega$, write
$P^\omega_x$ for the law of $Y$ started from $x\in\mathcal{C}_1$.
Moreover, define $Z$ to be the trace of $Y$ on $\mathcal{C}_2$, that
is, the
time change of $Y$ by the inverse of $A_t = \int_0^t 1_{( Y_s \in
\mathcal{C}
_2)} \,ds$. Specifically, writing ${\mathfrak a}_t=\inf\{s\dvtx A_s>t\}$
for the
right-continuous inverse of $A$, we set
\[
Z_t = Y_{{\mathfrak a}_t}, \qquad t\ge0.
\]
Note that unlike $Y$, the process $Z$ may perform long jumps by jumping
over the holes of $\mathcal{C}_2$. If $x \in\mathcal{C}_2(\omega)$
then we have
$Z_0=Y_0=x$, $P_x^\omega$-a.s., but otherwise \mbox{$Z_0=Y_{{\mathfrak a}_0}$}.

Given the percolation estimates of Section~\ref{percosect}, we can
follow \cite{ABDH}, Section~4, to establish the following theorems,
which correspond to Proposition 4.7(c) and Theorem 4.11 in \cite{ABDH}.
We remark that the second of the two results will be used in this paper
only for the proof of Theorem~\ref{teoEHIABDH}. Since it is the case
that, given Proposition~\ref{teoVNG-lem}, the proofs are a simple
modification of those in \cite{ABDH}, we omit them. For the statement
of the first result, we set
\[
\Psi(R,t) =\cases{ e^{ - R^2/t}, &\quad if $t> e^{-1} R$,
\cr
e^{ - R \log(R/t) }, &\quad if $t< e^{-1} R$.}
\]

\begin{propn}\label{teoldtqaiflfs} Write $\tau^Z_A = \inf\{t\dvtx Z_t
\notin A
\}$, $\tau^Y_A = \inf\{t\dvtx Y_t \notin A \}$. There exist constants
$\delta,c_i\in(0,\infty)$ and random variables $(R_x, x \in\mathbb
{L})$ with
%
\begin{equation}
\label{eqrxesprne} {\mathbb P}( R_x \ge n, x \in\mathcal{C}_1)
\le c_1 e^{-c_2 n^{\delta}},
\end{equation}
such that the following holds: for $x \in\mathcal{C}_1$, $t>0$ and $R
\ge R_x$,
\begin{eqnarray}
P_x^\omega\bigl( \tau^Z_{B_E(x,R)} < t
\bigr) &\le &c_3 \Psi(c_4 R,t),
\nonumber
\\
P_x^\omega\bigl( \tau^Y_{B_E(x,R)} < t
\bigr) &\le& c_3 \Psi(c_4 R,t).
\nonumber
\end{eqnarray}
\end{propn}

\begin{teo} There exist: constants $\delta,c_i\in(0,\infty)$; a set
$\Omega
_1\subset\Omega$ with ${\mathbb P}(\Omega_1)=1$; and random variables
$(S_x, x\in\mathbb{L})$ satisfying $S_x(\omega)<\infty$ for each
$\omega\in\Omega_1$ and $x\in\mathcal{C}_2(\omega)$, and
\[
{\mathbb P}( S_x \ge n, x \in\mathcal{C}_2) \le
c_1 e^{-c_2 n^{\delta}};
\]
such that the following statements hold.
\begin{longlist}[(a)]
\item[(a)] For $x,y\in\mathcal{C}_2(\omega)$ the transition density of $Z$,
as defined
by setting $q_t^Z(x,y):=P^{\omega}_x(Z_t=y)$, satisfies
\begin{eqnarray}
q^Z_t(x,y) &\le& c_3 t^{-d/2}\exp
\bigl(-c_4 |x-y|^2/t\bigr),\qquad t \ge|x-y|\vee
S_x,
\nonumber
\\
q^Z_t(x,y) &\ge& c_5 t^{-d/2}\exp
\bigl(-c_6 |x-y|^2/t\bigr),\qquad t \ge
|x-y|^{3/2}\vee S_x.
\nonumber
\end{eqnarray}
\item[(b)] Further, if $x \in\mathcal{C}_2(\omega)$, $t \ge S_x$ and
$B=B_2(x,2 \sqrt{t})$ then
\[
q^{Z,B}_t(x,y) \ge c_7t^{-d/2}\qquad
\mbox{for } y\in B_2(x,\sqrt{t}),
\]
where $B_2(x,R)$ is a ball in the (unweighted) graph $(\mathcal
{C}_2,\mathcal{O}_2)$, and $q^{Z,B}$ is the transition of $Z$ killed on
exiting $B$, that is, $q^{Z,B}_t(x,y):=P^{\omega}_x(Z_t=y,\tau_B^Z> t)$.
\end{longlist}
\end{teo}

Applying Proposition~\ref{teoldtqaiflfs}, we can establish the
following, which corresponds to \cite{ABDH}, Proposition 5.13(b). To
state the result, we introduce the rescaled process
$Y^{n}=(Y_t^n)_{t\geq0}$ by setting
\[
Y^n_t:=n^{-1} Y_{n^2t}.
\]

%
\begin{propn}[(Tightness)]\label{teotightABDH} Let $K,T,r>0$. For
${\mathbb P}$-a.e. $\omega$, the following is true: if $x_n\in D_n,
n\ge1$, $x\in B_E(0,K)$ are such that $x_n\to x$, then
%
\begin{eqnarray}
\label{eqtight1}  \lim_{R \to\infty}\limsup_{n\to\infty}
P_{nx_n}^\omega\Bigl( \sup_{s \le T} \bigl|
Y^{n}_s\bigr| > R \Bigr) &=& 0,
\\
\label{eqtight2}  \lim_{\delta\to0} \limsup_{n\to\infty}
P_{nx_n}^\omega\Bigl( \sup_{|s_1-s_2| \le\delta, s_i \le T } \bigl|
Y^{n}_{s_2}- Y^{n}_{s_1}\bigr| > r \Bigr)&=&
0.
\end{eqnarray}
In particular, for ${\mathbb P}$-a.e. $\omega$, if $x_n\in D_n, n\ge
1$, $x\in\overline{D}$ are such that $x_n\to x$, under~$P^\omega_{nx_n}$,
the family of processes $(Y_t^{n})_{t\geq0}, n\in\mathbb{N}$ is tight
in $\mathbb{D} ([0, \infty), \overline{D})$.
\end{propn}

\begin{pf} Since the statement is slightly different from \cite{ABDH},
Proposition 5.13, we sketch the proof. Note that since $x_n\to
x\in B_E(0,K)$, then by setting $M=K+1$ we have that $nx_n\in
B_E(0,nM)$ for all $n$ suitably large. Let $nR > \sup_nR_{nx_n}$. Then,
by Proposition~\ref{teoldtqaiflfs},
\[
P_{nx_n}^\omega\Bigl( \sup_{s \le T} \bigl|
Y^{n}_s\bigr| > R \Bigr) = P_{nx_n}^\omega
\bigl( \tau^Y_{B_E(0, nR)} < n^2T \bigr) \le
c_1 \Psi\bigl( c_2 nR, n^2T\bigr).
\]
Considering separately the cases $1/n < T/R$ and $1/n\ge T/R$, we
deduce that
\[
P_{nx_n}^\omega\Bigl( \sup_{s \le T} \bigl|
Y^{n}_s\bigr| > R \Bigr) \le c_3
e^{-c_4R^2/T } \vee e^{- R }.
\]
Since $\limsup_nR_{nx_n}/n\le\limsup_n\sup_{x\in
B_E(0,Mn)}R_{x}/n<\infty$, $\mathbb{P}$-a.s., due to the Borel--Cantelli
argument using (\ref{eqrxesprne}), we obtain (\ref{eqtight1}).

We next prove (\ref{eqtight2}). Write
\[
p(x,T,\delta, r) = P_{x}^\omega\Bigl( \sup
_{|s_1-s_2| \le\delta, s_i
\le T } | Y_{s_2}- Y_{s_1}| > r \Bigr),
\]
so that
\[
P_{nx_n}^\omega\Bigl( \sup_{|s_1-s_2| \le\delta, s_i \le T } \bigl|
Y^{n}_{s_2}- Y^{n}_{s_1}\bigr| > r \Bigr) =
p\bigl(nx_n, n^2T, n^2\delta,nr\bigr).
\]
Arguing similar to the proof of \cite{ABDH}, Proposition 5.13, we have
\[
p\bigl(nx_n, n^2T, n^2\delta, 2 n r\bigr)
\le c \exp\bigl(-c n T^{1/2} \bigr) + c (T/\delta) \exp\bigl( -c
r^2/\delta\bigr),
\]
provided
%
\begin{equation}
\label{z-ub2w2}\qquad T^{1/2} \ge n^{-1} R_x^{2/3},
\qquad\delta> n^{-1} r,\qquad r \ge n^{-1} \max
_{y \in B_E(nx_n, n^{3/2}T^{3/4})} R_y.
\end{equation}
Note that $B_E(nx_n, n^{3/2}T^{3/4})\subset B_E(0, nM+n^{3/2}T^{3/4})$
for large $n$.
If $T, r$ and $\delta$ are fixed, due to the Borel--Cantelli argument using
(\ref{eqrxesprne}), each of conditions in~(\ref{z-ub2w2}) holds when
$n$ is large enough. So, for $\mathbb{P}$-a.e. $\omega$,
\[
\limsup_{n\to\infty} p\bigl(nx_n, n^2T,
n^2\delta, 2 nr\bigr) \le c (T/\delta) \exp\bigl( -c r^2/
\delta\bigr)
\]
and (\ref{eqtight2}) follows.

Using (\ref{eqtight1}) and (\ref{eqtight2}), we have tightness for
$\{
P^{\omega}_{nx_n}\}$ by \cite{EK}, Corollary 3.7.4.
\end{pf}

We can further establish the following theorems, which correspond to
\cite{BiskP}, Lemma 5.6, Proposition 6.1 and \cite{ABDH}, Theorem 7.3.
We denote by\break $(q_t^Y(x,y))_{x,y\in\mathcal{C}_1,t>0}$ the heat kernel
associated with $Y$, that is, for $x,y\in\mathcal{C}_1$, $t>0$,
\[
q_t^Y(x,y):={P^{\omega}_x(Y_t=y)}.
\]
(We recall that the invariant measure of $Y$ is the uniform measure
$\nu
$ on $\mathcal{C}_1$.)

\begin{propn}\label{teobiskpreLem} There exist $c_1,c_2,c_3,\gamma\in
(0,\infty)$ (nonrandom) and random variables $(R_x, x\in\mathbb
{L})$ with
%
\begin{equation}
\label{eqdnwiw1334} \mathbb{P}(R_x\ge n, x\in{\mathcal C}_1)
\le\exp\bigl(-c_1n^{\gamma}\bigr),
\end{equation}
such that if $x,y\in{\mathcal C}_1$, then
\[
q_t^Y(x,y)\le c_2 t^{-d/2} \qquad
\mbox{for } t\ge\bigl(c_3\vee2d_1(x,y)\vee
R_x\bigr)^{1/4}.
\]
\end{propn}

\begin{pf} Note that the corresponding result for $Z$-process is given
in \cite{ABDH}, Corollary 4.3. We need to obtain similar result for
$Y$-process. First, note that
because we have
Proposition~\ref{teoVNG-lem}, the proof of \cite{ABDH}, Proposition 4.1
and Corollary 4.3 (with $\varepsilon=1/4$ for simplicity) goes
through once (4.7) in \cite{ABDH} is verified.
To check~\cite{ABDH}, (4.7),
we use \cite{BarDeu08}, Theorem 2.3, which can be proved almost
identically in our case. Note that in \cite{BarDeu08}, Theorem 2.3, the
metric $\tilde d$ is used, but thanks to Remark~\ref{teoremddtil}, we
can obtain the same estimates using the metric $d_1$. Finally, using
Cauchy--Schwarz, we obtain the desired inequality.
\end{pf}

For $G \subset\mathcal{C}_1$, we define $\partial^{\mathrm{out}}(G)$ to
be the
exterior boundary of $G$ in the graph $(\mathcal{C}_1, \mathcal
{O}_1)$, that
is, those vertices of $\mathcal{C}_1\setminus G$ that are connected to
$G$ by an edge in $\mathcal{O}_1$, and set $cl(G)= G \cup\partial
^{\mathrm{out}}(G)$. We say that a function $h$ is $Y$-harmonic in $A
\subset
\mathcal{C}
_1$ if $h$ is defined on $cl(A)$ and ${\mathcal L}_V h(x)=0$ for $x \in A$.

%
\begin{teo}[(Elliptic Harnack inequality)] \label{teoEHIABDH} There exist
random variables $(R'_x, x \in\mathbb{L})$ with
\[
{\mathbb P}\bigl( x \in\mathcal{C}_1, R'_x
\ge n\bigr) \le c e^{-c' n^\delta}
\]
and\vspace*{2pt} a constant $C_E$ such that if $x_0 \in\mathcal{C}_1$, $R \ge
R'_{x_0}$ and
$h\dvtx cl( B_1(x_0,R)) \to\mathbb{R}_+$ is $Y$-harmonic on $B_1=B_1(x_0,R)$,
then writing $B'_1= B_1(x_0,R/2)$,
\[
\sup_{B'_1} h \le C_E \inf_{B'_1}
h.
\]
\end{teo}

\begin{pf} Given Lemma~\ref{holelem}, the proof is almost identical
to that of \cite{ABDH}, Theorem 7.3.
\end{pf}

\subsection{Quenched invariance principle}\label{sec3.3}\label{qipsec}

To prove a quenched invariance principle for $Y$ (see Theorem~\ref{RCMresult} below), we will check the conditions
of Theorem~\ref{teosilvs}
one by one. We choose $\delta_n=c_*^2/n$, where $c_*$ is a constant
that will be chosen later. First, since $X^{n}$ is a continuous time
Markov chain with holding time at $x$ being an exponential random
variable of mean $\mu_x^{-1}$, it is conservative. Condition~(ii) in
Theorem~\ref{teosilvs} is a consequence\vspace*{2pt} of the quenched invariance
principle for the whole space (cf. \cite{ABDH}) and the fact that
$\mathcal{C}_1\subseteq\widetilde{\mathcal{C}}_1$, which is a consequence
of the uniqueness of the infinite percolation clusters in the two
settings. Since in the original papers quenched invariance principles
are uniformly stated in terms of the random walk started from the
origin, whereas we require such to hold from an arbitrary starting
point, we state\vspace*{2pt} the following generalization of existing results. We
will suppose that $\widetilde{P}_y^\omega$ refers to the quenched law of
the VSRW $\widetilde{Y}$ on $\widetilde{\mathcal{C}}_1$ started from
$y\in
\widetilde{\mathcal{C}}_1$, and $\widetilde{Y}^n$ refers to the rescaled
process defined by setting $\widetilde{Y}^n_t:=n^{-1}\widetilde
{Y}_{n^2t}$. The
proof of the result can be found in the \hyperref[sec5]{Appendix}.

\begin{teo}\label{qip} There exists a deterministic constant $c\in
(0,\infty)$
such that, for $\mathbb{P}$-a.e. $\omega$, the laws of the processes
$\widetilde{Y}^n$ under $\widetilde{P}_{nx_n}^\omega$, where $nx_n\in
\mathcal
{C}_1$ and \mbox{$x_n\rightarrow x$}, converge weakly to the laws of
$(B_{ct})_{t\geq0}$, where $(B_t)_{t\geq0}$ is standard Brownian
motion on $\mathbb{R}^d$ started from $x$.
\end{teo}

Concerning Assumption~\ref{teoassumpHC}(i), we will prove the
following: there exist $c_*, c_1\in(0,\infty)$ (nonrandom) and
$N_0(\omega)$ such that for all $n\ge N_0(\omega)$, $x_0,x \in
B_E(0,n^{1/2})\cap D_n$ and $c_*/n^{1/2}\le r'\le1$,
%
\begin{equation}
\label{eqimportant831} E^\omega_x \bigl(\tau_{B_E(x_0,r')\cap D_n}
\bigl(Y^n\bigr) \bigr)\le c_1r'^2.
\end{equation}
Applying Proposition~\ref{teobiskpreLem}, we have the following: for
all $x_0, x\in B_E(0, n^{3/2})\cap{\mathcal C}_1$ and $r\le n$, if
$c_2r^2\ge(c_3\vee2\sup_{z\in B_E(x_0,r)\cap{\mathcal
C}_1}d_1(x,z)\vee R_x)^{1/4}$, then
%
\begin{eqnarray}\label{eq13-1-1}
P^\omega_x \bigl(\tau^Y_{B_E(x_0,r)\cap{\mathbb L}}\ge
c_2r^2 \bigr)&\le& P^\omega_x
\bigl(Y_{c_2r^2}\in B_E(x_0,r) \bigr)\nonumber
\\
&=& \int
_{B_E(x_0,r)}q^Y_{c_2r^2}(x,z)\mu^0
(dz)
\\
&\le& \frac{c}{(c_2r^2)^{d/2}}\mu^0\bigl(B_E(x_0,r)
\bigr)\le\frac
{cc_4r^d}{c_2^{d/2} r^d}\le1/2,\nonumber
\end{eqnarray}
where\vspace*{1pt} we used $\mu^0(B_E(x_0,r))\le c_4r^d$ and we set
$c_2:=(2cc_4)^{2/d}\vee c_3^{1/4}$. Now, using Lemma~\ref{ballcomp},
there exists $N_1(\omega)\in{\mathbb N}$ that satisfies $\mathbb
{P}(N_1(\omega)\geq m)\leq c_1e^{-c_2m}$ such that $B_E(x,c_1R)\subset
B_1(x,R)$ for $x\in B_E(0,R)\cap{\mathcal C}_1$ and\break $R\ge N_1(\omega
)$. On the other hand, $|x-z|\le|x|+|x_0|+|x_0-z|\le2n^{3/2}+r$\break for
$z\in B_E(x_0,r)$, so taking $r\ge c_*n^{1/2}$ with $c_*$ large enough,
there\break exists $N_2(\omega)$  that satisfies
$\mathbb{P}(N_2(\omega)\geq m)\leq c_1e^{-c_2m}$ such that $c_2r^2\ge\break
(2\sup_{z\in B_E(x_0,r)\cap{\mathcal C}_1}d_1(x,z))^{1/4}$ holds for
$n\ge N_2(\omega)$. Next, by (\ref{eqdnwiw1334}),
\[
{\mathbb P} \Bigl( \sup_{x \in B_E(0,n^{3/2})\cap\mathcal{C}_1} R_x \ge
n \Bigr)\le
n^{3d/2} e^{-c_1 n^{\gamma}}.
\]
Summarizing, (\ref{eq13-1-1}) holds for all $x_0,x \in B_E(0,
n^{3/2})\cap{\mathcal C}_1$, $c_*n^{1/2}\le r\le n$ and $n\ge
N_0(\omega
):=N_2(\omega)\vee N_3(\omega)$, where $N_3(\omega):=\sup_{x \in
B_E(0,n^{3/2})\cap\mathcal{C}_1} R_x$. Moreover, the random variable
$N_0(\omega)$ is almost-surely finite; in fact, we have the following
tail bound for it, which will be useful in Example~\ref{boxsec} below:
%
\begin{equation}
\label{N0tail} \mathbb{P} \bigl(N_0(\omega)\geq m \bigr)\leq
c_1e^{-c_2m^\gamma}.
\end{equation}
Using the Markov property, we can inductively obtain
\[
P^\omega_x \bigl(\tau^Y_{B_E(x_0,r)\cap{\mathbb L}}\ge
kc_2r^2 \bigr)\le(1/2)^k\qquad\forall k\in
\mathbb{N}.
\]
So,
\begin{eqnarray*}
E^\omega_x \bigl(\tau^Y_{B_E(x_0,r)\cap{\mathbb L}} \bigr)
&\le& \sum_k (k+1)c_2r^2P^\omega_x
\bigl(kc_2r^2\le\tau^Y_{B_E(x_0,r)\cap
{\mathbb
L}}<
(k+1)c_2r^2 \bigr)
\\
&\le& 3c_2r^2.
\end{eqnarray*}
For $Y^n=n^{-1}Y_{n^2t}$, we therefore have: for $n\ge N_0(\omega)$,
$x_0\in B_E(0,n^{1/2})\cap D_n$, $x\in B_E(0,n^{1/2})\cap D_n$ and
$c_*/n^{1/2}\le r'\le1$,
\[
E^\omega_{nx} \bigl(\tau_{B_E(x_0,r')\cap D_n}
\bigl(Y^n\bigr) \bigr)\le\frac{1}{n^2}E^\omega_{nx}
\bigl(\tau^Y_{B_E(nx_0,r)\cap{\mathbb L}} \bigr)\le\frac{1}{n^2}
\cdot3c_2r^2=3c_2r'^2,
\]
where $r=nr'$. Thus, (\ref{eqimportant831}) holds, and so Assumption
\ref{teoassumpHC}(i) holds with $\delta_n=c_*^2/n$, $\beta=2$.

Regarding Assumption~\ref{teoassumpHC}(ii) we observe that, using
Theorem~\ref{teoEHIABDH}, the relevant condition can be obtained
similar to \cite{BarHam07}, Proposition~3.2. (Note that Proposition~3.2 in \cite{BarHam07} is a parabolic version, whereas we just need an
elliptic version.) Indeed, taking $(\log n)^{2/\delta}$ as $n$ in
Theorem~\ref{teoEHIABDH},
\[
{\mathbb P} \Bigl( \sup_{x \in B_E(0,cn^2)\cap\mathcal{C}_1} R'_x
\ge(\log n)^{2/\delta} \Bigr)\le c^dn^{2d}
e^{-c' (\log n)^{2}}\le c'/n^2.
\]
Thus, by the Borel--Cantelli lemma, there exists $N_1(\omega)\in
\mathbb{N}$
such that
\[
\sup_{x \in B_E(0,cn^2)\cap\mathcal{C}_1} R'_x \le(\log
n)^{2/\delta
}\qquad\forall n\ge N_1(\omega),
\]
so the elliptic Harnack inequality holds for $Y$-harmonic functions on
balls $B_E(x_0,R)$ with $x_0\in B_E(0,cn^2)$, $R\ge(\log n)^{2/\delta
}$. By scaling $Y^{n}(t)=n^{-1}Y_{n^2t}$, the elliptic Harnack
inequality holds uniformly for $Y^{n}$-harmonic functions on
$B_E(x_0,R)$ with $x_0\in B_E(0,cn)$, $R\ge(\log n)^{2/\delta}/n$.
Given the elliptic Harnack inequality, we can obtain the desired H\"
older continuity in a similar way as in the proof of Proposition 3.2 of
\cite{BarHam07}. Thus, setting $\delta_n:=c_*^2/n$, Assumption~\ref{teoassumpHC}(ii) holds for $R\ge\delta_n^{1/2}$, since $\delta
_n^{1/2}\ge(\log n)^{2/\delta}/n$.

Next, we remark that part (i) of Assumption~\ref{teoassumpHC-2} is
direct from Proposition~\ref{teotightABDH}, and part (ii) follows from
Proposition~\ref{teotightABDH} [especially (\ref{eqtight2}) implies
the condition].

The following proposition gives the appropriate convergence for the
sequence of measures $(m_n)_{n\geq1}$ defined by setting
$m_n:=n^{-d}\nu(n \cdot)$. (Recall that $\nu$ is the invariant measure
for $Y$, and so the measure $m_n$ is invariant for $Y^n$.)

\begin{propn}\label{vague} $\mathbb{P}$-a.s., the measures $(m_n)_{n\geq1}$
converge vaguely to $m$, a~deterministic multiple of Lebesgue measure
on $\overline{D}$.
\end{propn}

\begin{pf} First note that if $Q \subset\overline{D}$ is a cube of side
length $\lambda$, then applying Lemma~\ref{meascomp} in a Borel--Cantelli argument yields that, $\mathbb{P}$-a.s.,
%
\begin{equation}
\label{n1n2} \frac{\tilde{\nu}(n Q)-\nu(nQ)}{n^d}\rightarrow0.
\end{equation}
Next, consider a rectangle of the form $R=[0,\lambda_1]\times\cdots
\times
[0,\lambda_d]$. Since the full $\mathbb{Z}^d$ model is ergodic under
coordinate shifts, a simple application of a multidimensional ergodic
theorem yields that, $\mathbb{P}$-a.s.,
\[
\frac{\tilde{\nu}(nR)}{n^d}\rightarrow c_1 \prod
_{i=1}^d\lambda_i,
\]
where $c_1:=\mathbb{P}(0\in\widetilde{\mathcal{C}}_1)\in(0,1]$. An
inclusion--exclusion argument allows one to extend this result to any
rectangle of the form $[x_1,x_1+\lambda_1]\times\cdots\times
[x_d,x_d+\lambda_d]$, where $x_i\geq0$ for $i=1,\ldots,d$. Clearly, the
particular orthant is not important, so the result can be further
extended to cover any rectangle $R\subset\overline{D}$, and, in particular,
we have that, $\mathbb{P}$-a.s.,
\[
\frac{\tilde{\nu}(n Q)}{n^d}\rightarrow c_1\lambda^d.
\]
Applying (\ref{n1n2}), we obtain that the above limit still holds when
$\tilde{\nu}$ is replaced by~$\nu$. The proposition follows.
\end{pf}

Finally, in order to verify condition~(iii) in Theorem~\ref{teosilvs},
we first give a lemma.

\begin{lem}\label{teoLLN-lem} Let $\{\eta_i\}_{i\geq1}$ be independent and
identically distributed with $E|\eta_1|<\infty$. Suppose $\{a_k^n\}
_{k=1}^n$ is a sequence of real numbers with $|a_k^n|\le M$ for all
$k,n$ (here $M>0$ is some fixed constant) such that the following two
limits exist:
\[
a:=\lim_{n\to\infty}\frac{1}n \sum
_{k=1}^na_k^n,\qquad\lim
_{n\to
\infty
}\frac{1}n \sum_{k=1}^n\bigl|a_k^n\bigr|.
\]
It then holds that $\frac{1}n \sum_{k=1}^na_k^n\eta_k$ converges to
$aE[\eta_1]$ almost surely as $n\to\infty$.
\end{lem}

\begin{pf} This can be proved similarly to Etemadi's proof of the
strong law of large numbers (see \cite{Dur}).
\end{pf}

\begin{propn}\label{teoiiiinthsil}
If $C_0:=\mathbb{E}(\mu_e)<\infty$, then
$\mathbb{P}$-a.s.,
%
\begin{equation}
\qquad \widetilde{\mathcal{E}}(f,f)\le\limsup_{n\to\infty}\mathcal
{E}^{(n)}(f,f)\leq C_0\int_{\overline{D}}\bigl|
\nabla f(x)\bigr|^2\,dx \qquad\forall f\in C_c^2(
\overline{D}),\label
{eqwqdf211}
\end{equation}
where $C_c^2(\overline{D})$ is the space of compactly supported
functions on
$\overline{D}$ with a continuous second derivative, and
\[
\mathcal{E}^{(n)}(f,f):=\frac{n^{2-d}}2\mathop{\sum
_{x,y\in
{\mathcal
{C}_1}\dvtx}}_{\{x,y\}\in\mathcal{O}_1}\bigl(f(x/n)-f(y/n)\bigr)^2
\mu_{x,y}
\]
is the Dirichlet form on $L^2(D_n; m_n)$ corresponding to $Y^n$. In
particular, condition~\textup{(iii)} in Theorem~\ref{teosilvs} holds.
\end{propn}

\begin{pf} The first inequality of (\ref{eqwqdf211}) is standard. Indeed,
\begin{eqnarray*}
\widetilde{\mathcal{E}}(f,f)&=& \sup_{t>0}\frac{1}t
(f-P_tf,f)=\sup_{t>0}\liminf_{k\to
\infty}
\frac{1}t \bigl(f-P^{n_{j(k)}}_tf,f\bigr)
\\
&\le& \liminf_{k\to\infty}\sup_{t>0}
\frac{1}t \bigl(f-P^{n_{j(k)}}_tf,f\bigr)=\liminf
_{k\to\infty}\mathcal{E}^{(n_{j(k)})}(f,f),
\end{eqnarray*}
where the first inner product is in $L^2(\overline{D}; m)$, and the
other two are in $L^2(D_{n_{j(k)}}; m_{n_{j(k)}})$. Moreover, to
establish the second
inequality of (\ref{eqwqdf211}),
we apply the local uniform convergence
of $P^{n_{j(k)}}_tf$ to $P_tf$ (cf. the proof of Proposition~\ref{clttight}) and the vague convergence of $m_{n_{j(k)}}$ to $m$ (Lemma
\ref{vague}). We now prove the second inequality. Suppose $\operatorname{Supp}f\subset B_E(0,M)\cap\overline{D}$ for some $M>0$, then
\[
\mathcal{E}^{(n)}(f,f)=\frac{n^{2-d}}2\sum
_{(x,y)\in
H_{n,M}}\bigl(f(x/n)-f(y/n)\bigr)^2
\mu_{x,y},
\]
where $H_{n,M}:=\{(x,y)\dvtx x,y\in B_E(0,nM)\cap\mathcal{C}_1, \{x,y\}
\in
\mathcal{O}_1\}$. Clearly, this
quantity increases when $H_{n,M}$ is replaced by
$H'_{n,M}:=\{(x,y)\dvtx x\in
B_E(0,nM)\cap\mathbb{L}, \{x,y\}\in E_{\mathbb{Z}^d}\}$.
Note that $\# H'_{n,M}\sim c_1 (nM)^d$. Moreover, set
$a_{(x,y)}^n:=\break n^2(f(x/n)- f(y/n))^2$ and $\eta_{(x,y)}:=\mu_{x,y}$.
Then, since $f\in C_c^2$, $0\le a_{(x,y)}\le M'$ for some
$M'>0$, and further,
\[
\lim_{n\to\infty}\bigl(2c_1(nM)^d
\bigr)^{-1}\sum_{(x,y)\in
H'_{n,M}}a_{(x,y)}^n=c_1^{-1}M^{-d}
\int_{\overline{D}}\bigl|\nabla f(x)\bigr|^2\,dx
\]
by applying Lemma~\ref{teoLLN-lem}, we obtain that, $\mathbb{P}$-a.s.,
\[
\lim_{n\to\infty}n^{-d}\sum_{(x,y)\in H'_{n,M}}a_{(x,y)}^n
\eta_{(x,y)}=2C_0\int_{\overline{D}}\bigl|\nabla
f(x)\bigr|^2\,dx.
\]
The result at (\ref{eqwqdf211}) follows.
\end{pf}

Putting together the above results, we conclude the following.

\begin{teo}\label{RCMresult} For the random conductance model on
$\mathbb{L}$
with independent and identically distributed conductances $(\mu
_e)_{e\in E_{\mathbb{L}}}$ satisfying (\ref{mucond1}), (\ref{mucond2})
and (\ref{mucond3}), there exists a deterministic constant $c\in
(0,\infty)$ such that, for $ \mathbb{P}$-a.e. $\omega$, the laws of
the processes ${Y}^n$ under ${P}_{nx_n}^\omega$, where $nx_n\in
\mathcal
{C}_1$ and $x_n\rightarrow x$, converge weakly to the laws of
$\{X_{ct}; t\geq0\}$,
where $\{X_t; t\geq0\}$ is the reflecting Brownian motion on
$\overline{D}$ started from $x$.
\end{teo}

\begin{rem}\label{genrem} (i) The diffusion constant $c$ is the same for
the model restricted to $\mathbb{L}$ as for the full $\mathbb{Z}^d$
model.

%
\begin{figure}[t]

\includegraphics{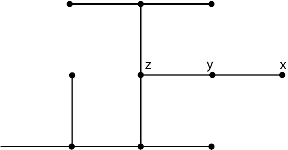}

\caption{An example trap structure.}\label{trapfig}
\end{figure}

(ii) When the conductance is not bounded from below, we cannot apply
our theorem because Assumption~\ref{teoassumpHC}(i) does not hold in
general, and we do not know how to obtain the quenched invariance
principle without this. Indeed, consider the realization of edge
weights shown in Figure~\ref{trapfig}, where the conductance on $\{
x,y\}
$ is $1$ and it is $O(n^{-\alpha})$ on $\{y,z\}$ where $\alpha>2$. One
can easily compute that $E^\omega_x\tau_{B_E(x,2)}(Y)\ge c_1
n^{\alpha}
\gg n^2$. Let $p_0=\mathbb{P} (\mu_e=0 )$ and $p_1=\mathbb
{P} (\mu_e=1 )$. Then the probability that such a trap
configuration appears is $p_0^{4d-3}p_1\mathbb{P} (0<\mu_e\le
n^{-\alpha} )=c_2\mathbb{P} (0<\mu_e\le n^{-\alpha
} )$.
Now let $\Omega_n:=\{\exists x_n\in B_E(0,n/2)$ such that $E^\omega_{x_n}\tau_{B_E(x_n,2)}(Y)\ge c_3n^{\alpha}\}$. If we have
$\mathbb{P} (0<\mu_e\le x )\ge c_4x^{d/\alpha}$ for small
$x>0$, then
$\mathbb P(\Omega_n)\ge1-(1-c_4n^{-d})^{n^d}\ge1-e^{-c_5}$ for large
$n$. In particular, $\limsup_n \Omega_n$ occurs with positive
probability. Set $X^n_t:=n^{-1}Y_{n^2t}$. Then, for $\omega\in\Omega
_n$, we have
\begin{eqnarray*}
E^\omega_{n^{-1}x_n} \bigl(\tau_{B_E(0,1)\cap D_n}
\bigl(X^n\bigr) \bigr)&=&E^\omega_{x_n} \bigl(
\tau_{B_E(0,n)\cap D_0}(Y_{n^2\cdot}) \bigr)
\\
&\ge& E^\omega
_{x_n} \bigl(\tau_{B_E(x_n,2)}(Y_{n^2\cdot}) \bigr)\ge
c_3n^{\alpha-2}.
\end{eqnarray*}
Since $\Omega_n$ occurs infinitely often with positive probability,
Assumption~\ref{teoassumpHC}(i)
does not hold for any choice of $\beta>0$ (by choosing $x_0=0, r=1$).

(iii) There is another natural continuous time Markov chain on
$\mathcal
{C}_1$, namely with transition probability $P(x,y)=\mu_{xy}/\mu_x$ and
the holding time being the exponential distribution with mean one for
each point. [Such a Markov chain is sometimes called a constant speed
random walk (CSRW).] It is a time change of the VSRW; the corresponding
Dirichlet form is $(\mathcal{E},L^2(\mathcal{C}_1; \mu))$, and the
corresponding discrete Laplace operator is ${\mathcal L}_C f(x)=\frac
{1}{\mu_x}\sum_y(f(y)-f(x))\mu_{xy}$. For the\vspace*{1.5pt} whole space case, one can
deduce the quenched invariance principle of CSRW from that of VSRW by
an ergodic theorem. (See \cite{ABDH}, Section~6.2 and \cite{BarDeu08},
Section~5. Note that the limiting process degenerates if
${\mathbb E}\mu_e=\infty$.) Since our state space $\mathbb L$ features
a lack of translation invariance, we cannot use the ergodic theorem. So
far, we do not know how to circumvent this issue to prove the quenched
invariance principle for general CSRW on $\mathbb L$. (However, we do
note that for the case of random walk on a supercritical percolation
cluster, the CSRW and VSRW behave similarly, and the quenched
invariance principle for the CSRW can be proved in a similar way as for
the VSRW case. Moreover, the quenched invariance principle for the
discrete time simple random walk on $\mathcal{C}_1$ follows easily from
that for the CSRW.)

(iv) To extend Theorem~\ref{RCMresult} to apply to more general
domains, it will be enough to check the percolation estimates from
which we deduced Assumptions~\ref{teoassumpHC} and~\ref{teoassumpHC-2} in these settings. While we believe doing so should be
possible, at least under certain smoothness assumptions on the domain
boundary, we do not feel the article would benefit significantly by the
increased technical complication of pursuing such results, and
consequently omit to do so here. Instead, we restrict our discussion of
more general domains to the case of uniformly elliptic conductances,
where the relevant estimates are straightforward to check (see
Section~\ref{uesec} below). Similarly, given suitable full-space
quenched invariance principles and percolation estimates (namely the
estimates given in Lemmas~\ref{holelem}--\ref{ballcomp}), our results
should readily adapt to percolation models on other lattices.

(v) Given the various estimates we have established so far, it is
possible to extend the quenched invariance principle of Theorem~\ref{RCMresult} to a local limit theorem, that is, a result describing the
uniform convergence of transition densities. More specifically, the
additional ingredient needed for this is an equi-continuity result for
the rescaled transition densities on $\mathcal{C}_1$, which can be
obtained by applying an argument similar to that used to deduce
Assumption~\ref{teoassumpHC}(ii), together with the heat kernel upper
bound estimate
of Proposition~\ref{teobiskpreLem}. Since the proof of such a result
is relatively standard (cf. \cite{BarHam07,CroHam08}), we will only
write out the details in the compact box case (see Section~\ref{boxsec}
below), where convergence of transition densities is also useful for
establishing convergence of mixing times.
\end{rem}

\section{Other examples}\label{sec4}\label{exsec}

\subsection{Random conductance model in a box}\label{sec4.1}\label{boxsec}

The purpose of this section is to explain how to adapt the results of
Section~\ref{percsec} to the compact space case. For $d\geq2$ fixed,
set $\mathbb{B}(n):=[-n,n]^d\cap\mathbb{Z}^d$, let ${E}_{\mathbb
{B}(n)}=\{e=\{x,y\}\dvtx x,y\in\mathbb{B}(n), |x-y|=1\}$ be the set of
nearest neighbor bonds, and suppose $\mu=(\mu_e)_{e\in E_{\mathbb
{Z}^d}}$ is a collection of independent random variables satisfying the
assumptions made on the weights in Section~\ref{percsec}, that
is, (\ref{mucond1}), (\ref{mucond2}) and (\ref{mucond3}). For each $n$
and each
realization of $\mu$, let $\mathcal{C}_1(n)$ be the largest component
of $\mathbb{B}(n)$ that is connected by edges satisfying $\mu_e>0$, and
let $Y^n$ be the VSRW on $\mathcal{C}_1(n)$. We will write
$P_{n,x}^\omega$ for the quenched law of $Y^n$ started from $x\in
\mathcal{C}_1(n)$. The aim of this section is to show, via another
application of Theorem~\ref{teosilvs}, that $X^n=(X_t^n)_{t\geq0}$,
defined by setting
\[
X^n_t:=n^{-1}Y^n_{n^2t},
\]
converges as $n\rightarrow\infty$ to reflecting Brownian motion on
$D=[-1,1]^d$, for almost-every realization of the random environment
$\mu$. We observe that, in the case of uniformly elliptic random
conductances, this result was recently established using an alternative
argument in \cite{BuckBen}.
Note also that, by applying a result from \cite{CrHK}, the above
functional scaling limit readily yields the corresponding convergence
of mixing times (see Corollary~\ref{mtcor} below for a precise statement).

To prove the results described in the previous paragraph, we start by
considering a decomposition of $\mathbb{B}(n)$. In particular, fix
$\varepsilon\in(0,1)$ and for $i=(i_1,\ldots,i_d)\in\{-1,1\}^d$, let
$\mathbb{B}_i^\varepsilon(n)$ be the cube of side-length $\lfloor
n(1+\varepsilon)\rfloor$ which has a corner at $ni$ and contains 0.
Within each of the $2^d$ sets of the form $\mathbb{B}^\varepsilon
_i(n)$, the random walk on $\mathcal{C}_1(n)$ reflects only at the
faces of the box adjacent to the single corner vertex~$ni$. As a
consequence, we will be able to transfer a number of key estimates to
the current framework from the unbounded case considered in
Section~\ref{percsec}---note that the reason for taking $\varepsilon>0$ is so that
the boxes $\mathbb{B}_i^\varepsilon(n)$ overlap, which will allow us to
``patch'' together results proved for different parts of the box. For
the purpose of transferring results from Section~\ref{percsec}, the
following lemma will be useful. Its proof can be found in the \hyperref[sec5]{Appendix}.

\begin{lem}\label{boxclust} There exist constants $c_1,c_2$ such that if
$Q_1, Q_2\subseteq\mathbb{Z}_+^d$ are the cubes of side length
$\lfloor n(1+\varepsilon)\rfloor$, $2n$ containing 0, respectively,
$\mathcal{C}^+(Q_2)$ is the largest connected component of the graph
$(Q_2,\mathcal{O}_1)$, and $\mathcal{C}_1$ is the unique infinite
component of $(\mathbb{Z}_+^d,\mathcal{O}_1)$, then
\[
\mathbb{P} \bigl(\mathcal{C}^+(Q_2)\cap Q_1\neq
\mathcal{C}_1\cap Q_1 \bigr)\leq c_1e^{-c_2n}.
\]
\end{lem}

%
\begin{figure}[t]

\includegraphics{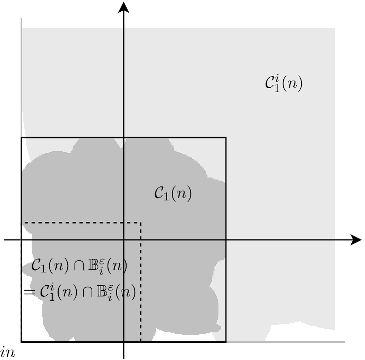}

\caption{Schematic illustration of how, for large $n$, the largest
cluster $\mathcal{C}_1(n)$ (dark-grey) in $\mathbb{B}(n)$ (square with
solid boundary) and the largest cluster $\mathcal{C}_1^i(n)$
(light-grey) in the quarter-plane with corner $in$ and containing 0
agree within $\mathbb{B}_i^\varepsilon(n)$ (square with dashed
boundary).}\label{boxfig}
\end{figure}

In particular, if $\mathcal{C}_1^i(n)$ is the unique infinite
percolation cluster on the copy of~$\mathbb{Z}_+^d$ that has corner
vertex $ni$ and contains 0, then the above result implies that with
probability at least $1-c_1e^{-c_2n}$ (or, by Borel--Cantelli,
almost-surely for large $n$) we have that
%
\begin{equation}
\label{inclusion} \mathcal{C}_1(n)\cap\mathbb{B}_i^\varepsilon(n)=
\mathcal{C}_1^i(n)\cap\mathbb{B}^\varepsilon_i(n)
\subseteq\widetilde{\mathcal{C}}_1\qquad\forall i\in\{-1,1
\}^d,
\end{equation}
where the inclusion is a consequence of the uniqueness of the infinite
percolation clusters in question. (This result is summarized by
Figure~\ref{boxfig}.)

We now check the conditions listed at the end of Section~\ref{sect2}
one by one. Since in light of (\ref{inclusion}), most of these are
straightforward adaptations of the arguments given in Section~\ref
{percsec}, we will only provide a brief description of how to do this.
First, as was the case in the $\mathbb{L}$ setting, since $X^{n}$ is a
continuous time Markov chain with holding time at $x$ being $\exp(\mu
_x)$, it is conservative. Second, given (\ref{inclusion}),
condition~(ii) in Theorem~\ref{teosilvs} is a consequence of the
quenched invariance principle for the whole space stated as Theorem
\ref
{qip}. Moreover, since $\mathcal{C}_1(n)$ agrees with $\widetilde
{\mathcal
{C}}_1$ up to a distance $c\log n$ of the boundary, at least for large
$n$ (see \cite{BenMos}, Proposition 1.2), by applying the full
$\mathbb
{Z}^d$ version of Proposition~\ref{vague}, we have that the measures~$m_n$, defined analogously to the previous section, $\mathbb
{P}$-a.s. converge weakly to (a suitably rescaled version of) Lebesgue
measure on $[-1,1]^d$. Similarly, the Dirichlet form comparison of
(iii) can be obtained by following the same argument used to prove the
corresponding result in Section~\ref{percsec}---Proposition~\ref{teoiiiinthsil}. Applying (\ref{inclusion}), we are also able to
deduce the following tightness result, which is analogous to
Proposition~\ref{teotightABDH}, and from which Assumption~\ref{teoassumpHC-2} is readily obtained.

\begin{propn} For ${\mathbb P}$-a.e. $\omega$, if $x_n\in n^{-1}\mathcal
{C}_1(n)$, $n\geq1$ is such that $x_n\to x\in\overline{D}$, then under
$P^\omega_{n,nx_n}$,
the family of processes $(X_t^{n})_{t\geq0}, n\in\mathbb{N}$ is  tight
in $\mathbb{D} ([0, \infty), [0, 1]^d)$, and any convergent subsequence
has limit in\break  $C ([0, \infty), [0, 1]^d)$.
\end{propn}

\begin{pf} Note that in the bounded case the limit corresponding to
(\ref{eqtight1}) is immediate, and hence it will suffice to check the
limit corresponding to (\ref{eqtight2}). To do this, the same argument
can be applied, so long as one can check the following: for any $r\in
(0,\varepsilon)$, there exist $c_i$ and random variables $(R_x^n, x\in
\mathbb{B}(n),n\geq1)$ with
\[
{\mathbb P}\bigl( R_x^n \ge r n, x \in
\mathcal{C}_1(n)\bigr) \le c_1 e^{-c_2
n^{\delta}},
\]
such that if $x \in\mathcal{C}_1(n)$, $t>0$ and $R \ge R_x^n$, then
%
\begin{equation}
\label{rbound}P_{n,x}^\omega\bigl( \tau^{Y^n}_{B_E(x,R)}
< t \bigr) \le c_3 \Psi(c_4 R,t).
\end{equation}
For this purpose, if $x\in\mathbb{B}^0_i(n)$, set $\widetilde{R}_x^n$ to
be equal to $R_x^{n,i}$, the quantity defined in Proposition~\ref{teoldtqaiflfs} with $\mathcal{C}_1$ replaced by $\mathcal{C}_1^i(n)$.
If $\widetilde{R}_x^n\leq\varepsilon n$ and the part of $\mathcal{C}_1(n)$
contained in $B_E(x,\varepsilon n)$ is identical to the part of
$\mathcal{C}^i_1(n)$ contained in this set, then set $R_x^n=\widetilde
{R}_x^n$. Else, set ${R}_x^n=3n$. The required exponential decay for
the distributional tail of $R_x^n$ then follows from Proposition~\ref{teoldtqaiflfs} and (\ref{inclusion}). Moreover, since the probability
on the left-hand side of (\ref{rbound}) is 0 for $R\geq3n$, the bound
at (\ref{rbound}) follows.
\end{pf}

It remains to check Assumption~\ref{teoassumpHC}. For part (i), we
simply note that the combination of (\ref{eqimportant831}) [or more
precisely, the exponential tail bound for $N_0$ that appears as (\ref
{N0tail})] and (\ref{inclusion}) in a standard Borel--Cantelli argument
implies the following: there exist $c_*, c_1>0$ (nonrandom) and
$N_0(\omega)$ such that if $n\geq N_0(\omega)$, then, for each $x_0,x
\in n^{-1}\mathcal{C}_1(n)$ and $c_*/n^{1/2}\le r'\le1$,
\[
E^\omega_{n,x} \bigl(\tau_{B_E(x_0,r')}
\bigl(X^{n}\bigr) \bigr)\le c_1r'^2
\]
as desired. For part (ii) of this assumption, first note that we can
obtain the elliptic Harnack inequality uniformly for $X^n$-harmonic
functions on $B_E(x_0,R)$, where $x_0\in n^{-1}\mathcal{C}_1(n)$ and
$(\log n)^{2/\delta}/n\leq R\leq1$ for large $n$. [This can be
proved similarly as before, namely when $x_0\in\mathbb{B}^0_i(n)$,
Theorem~\ref{teoEHIABDH} can be applied by replacing $\mathcal{C}_1$
by $\mathcal{C}_1^i(n)$ due to (\ref{inclusion}).] Given the elliptic
Harnack inequality, we can obtain H\"older continuity in a similar way
as in the proof of \cite{BarHam07}, Proposition 3.2, for example.
Hence, we have established the following.

\begin{teo}\label{boxresult} There exists a constant $c\in(0,\infty)$ such
that, for $\mathbb{P}$-a.e. $\omega$, the process $X^n$ under
${P}_{n,nx_n}^\omega$, where $nx_n\in\mathcal{C}_1$ and
$x_n\rightarrow
x\in[-1,1]^d$, converges in distribution to $(B_{ct})_{t\geq0}$, where
$(B_{t})_{t\geq0}$ is
the reflecting Brownian motion on $[-1,1]^d$ started from $x$.
\end{teo}

Next, for $p\in[1,\infty]$, define the $L^p$ mixing time of the VSRW on
$\mathcal{C}_1(n)$ to be
%
\begin{eqnarray}\label{mixdef}
&& t_{\mathrm{mix}}^p\bigl(\mathcal{C}_1(n)
\bigr)
\nonumber\\[-8pt]\\[-10pt]
&&\qquad:=\inf\biggl\{t>0\dvtx \sup_{x\in\mathcal{C}_1(n)} \biggl(\int
_{\mathcal{C}_1(n)}\bigl\llvert q^n_t(x,y)-1\bigr
\rrvert^p\pi^n(dy) \biggr)^{1/p}<
\frac{1}{4} \biggr\},\nonumber
\end{eqnarray}
where we denote by $q^n$ the transition density of the VSRW with
respect to its (unique) invariant probability measure $\pi^n$. The
above result then has the following corollary. Note that in the
percolation setting, the obvious adaptation of this result to discrete
time gives a refinement of the first statement of \cite{BenMos},
Theorem 1.1.

\begin{cor}\label{mtcor}
Fix $p\in[1,\infty]$. For $\mathbb
{P}$-a.e. $\omega
$, we have that
\[
n^{-2}t_{\mathrm{mix}}^p\bigl(\mathcal{C}_1(n)
\bigr)\rightarrow c^{-1}t_{\mathrm{mix}}^p
\bigl([-1,1]^d\bigr),
\]
where $c$ is the constant of Theorem~\ref{boxresult}, and $t_{\mathrm
{mix}}^p([-1,1]^d)$ is the mixing time of reflecting Brownian motion on
$[-1,1]^d$ [defined analogously to (\ref{mixdef})].
\end{cor}

\begin{pf} First note that a simple rescaling yields that, $\mathbb
{P}$-a.s., $\pi^n$ converges weakly to a rescaled version of Lebesgue
measure on $[-1,1]^d$. The $\mathbb{P}$-a.s. Hausdorff convergence of
$n^{-1}\mathcal{C}_1(n)$ (equipped with Euclidean distance) to
$[-1,1]^d$ is a straightforward consequence of this. To establish the
corollary by applying \cite{CrHK}, Theorem 1.4, it will thus be enough
to extend the weak convergence result of Theorem~\ref{boxresult} to a
uniform convergence of transition densities (so as to satisfy \cite
{CrHK}, Assumption 1). According to \cite{CrHK}, Proposition 2.4 (cf.
\cite{CroHam08}, Theorem 15) and the quenched invariance principle
mentioned above, it is enough to show \cite{CrHK}, (2.11), namely, for
any $0<a<b<\infty$,
%
\begin{equation}
\label{tightcond} \qquad\lim_{\delta\rightarrow0}\limsup_{n\rightarrow\infty}
\mathop{\sup_{x,y,z\in n^{-1}\mathcal{C}_1(n)\dvtx
}}_{d_{E}(ny,nz)\leq n\delta}\sup_{t
\in[a,b]}
\bigl\llvert q_{n^2t}^{n}(nx,ny)-q_{n^2t}^{n}(nx,nz)
\bigr\rrvert=0.
\end{equation}
To prove this, first we have the following
H\"older continuity, which can be checked similarly to Assumption~\ref{teoassumpHC}(ii):
%
\begin{eqnarray}\label{wnibfs}
\bigl\llvert q_{n^2t}^{n}(nx,ny)-q_{n^2t}^{n}(nx,nz)
\bigr\rrvert\leq c_1|y-z|^\gamma\bigl\|q_{n^2t}^{n}(nx,
\cdot)\bigr\|_\infty
\nonumber\\[-8pt]\\[-8pt]
\eqntext{\forall x,y,z\in n^{-1}
\mathcal{C}_1(n).}
\end{eqnarray}
For $0<a\le t<1$, say, a compact version of Proposition~\ref{teobiskpreLem} and scaling gives that
\[
\bigl\|q_{n^2t}^{n}(nx,\cdot)\bigr\|_\infty\le
c_3\bigl|\mathcal{C}_1(n)\bigr|\bigl(n^2t
\bigr)^{-d/2} \leq c_3a^{-d/2}
\]
for large $n$. For $t\ge1$, Cauchy--Schwarz and monotonicity of
$q_{n^2t}^{n}(nx,nx)$
implies $\|q_{n^2t}^{n}(nx,\cdot)\|_\infty\le c_4$. In particular,
%
\begin{equation}
\label{eqHKondg} \bigl\|q_{n^2t}^{n}(nx,\cdot)\bigr\|_\infty\le
c_2(a)
\end{equation}
uniformly in $x\in D_n$, $t\geq a$, for large $n$, $\mathbb{P}$-a.s.
Thus, for $t\in[a,b]$, the right-hand side of (\ref{wnibfs}) is
bounded from above by $c_3(a) |y-z|^\gamma$. Taking $n\to\infty$ and
then $\delta\to0$, we obtain (\ref{tightcond}).
\end{pf}

Finally, as a corollary of the heat kernel continuity derived in the
proof of the previous result, we obtain the following local central
limit theorem. We let $g_n\dvtx [-1,1]^d\to{\mathcal C}_1(n)$ be such that
$g_n(x)$ is a closest point in ${\mathcal C}_1(n)$ to $nx$ in the
$|\cdot|_\infty$-norm. (If there is more than one such point, we choose
one arbitrarily.)

\begin{cor} Let $q_t(\cdot,\cdot)$ be the heat kernel of the reflecting
Brownian motion on $[-1,1]^d$.
For $\mathbb{P}$-a.e. $\omega$ and for any $0<a<b<\infty$, we have that
%
\begin{equation}
\label{eqLCLTE} \lim_{n\to\infty}\sup_{x,y\in[-1,1]^d}\sup
_{t\in[a,b]} \bigl|q_{nt}^n\bigl(g_n(x),g_n(y)
\bigr)-q_{ct}(x,y)\bigr|=0,
\end{equation}
where $c$ is the constant of Theorem~\ref{boxresult}.
\end{cor}

\begin{pf} Given the above results, the proof is standard. By (\ref
{wnibfs}) and (\ref{eqHKondg}) and the Ascoli--Arzel\`a theorem (along
with the Hausdorff convergence of $n^{-1}\mathcal{C}_1(n)$ to
$[-1,1]^d$), there exists a subsequence of $q_{nt}^n(g_n(\cdot
),g_n(\cdot))$ that\vspace*{1pt} converges uniformly to a jointly continuous
function on $[a,b]\times[-1,1]^d\times[-1,1]^d$. Using Theorem~\ref{boxresult}, it can be checked that this function is the heat kernel of
the limiting process. Since the limiting process is unique, we have the
convergence of the full sequence of $q_{nt}^n(g_n(\cdot),g_n(\cdot))$.
The uniform convergence in (\ref{eqLCLTE}) is then another consequence
of (\ref{wnibfs}) and (\ref{eqHKondg}).
\end{pf}

\subsection{Uniformly elliptic random conductances in uniform domains}\label{sec4.2}\label{uesec}

When the conductances are uniform elliptic, that is, bounded from above
and below by fixed positive constants, we can obtain quenched
invariance principles for a much wider class of domains than those
considered in the examples presented so far. In particular, let $D$ be
a uniform domain in $\mathbb{R}^d$, $d\geq2$.
For each $n\ge1$, let $\widehat D_n$ be the largest connected component
of $nD\cap{\mathbb Z}^d$, and set ${E}_{\widehat D_n}=\{e=\{x,y\}\dvtx
x,y\in\widehat D_n, |x-y|=1\}$. Suppose $\mu=(\mu_e)_{e\in E_{\mathbb
{Z}^d}}$ is a collection of independent random variables such that
\[
\mathbb P(C_1\le\mu_e\le C_2)=1\qquad
\forall e\in{E}_{\mathbb{Z}^d},
\]
where $C_1,C_2$ are nonrandom positive constants. Let $Y^n$ be either
VSRW or CSRW on $\widehat D_n$. Moreover, set $D_n:=n^{-1}\widehat D_n$ and
define $X^n_t:=n^{-1}Y^n_{n^2t}$. It is then the case that
$X^n=(X_t^n)_{t\geq0}$ converges as $n\rightarrow\infty$ to a
(constant time change of) reflecting Brownian motion on $D$, for
$\mathbb{P}$-almost-every realization of the random environment~$\mu$.

Since checking the details for this case is much simpler than for
previous settings, we will not provide a full proof of the result
described in the previous \mbox{paragraph}, but merely indicate how to
establish the key estimates. Indeed, in the uniformly elliptic case, we have
\begin{eqnarray*}
c_1R^d &\le& n^{-d}\bigl\llvert
B_E(x,R)\cap D_n\bigr\rrvert\le c_2R^d,
\\
c_3R^d &\le& n^{-d}\mu
\bigl(B_E(x,nR)\cap\widehat D_n\bigr)\le
c_4R^d
\end{eqnarray*}
for all large $n$ and all $n^{-1}\le R\le\operatorname{diam} (D)$, $x\in D$,
where $c_i$ are nonrandom positive constants. Furthermore, the weak
Poincar\'{e} inequality (\ref{poincare}) holds both for the counting
measure and $n^{-d}\mu$ uniformly for $n^{-1}\le R\le\operatorname{diam}(D)$,
in the sense that the constants do not depend on $n$. Given these, we
can apply \cite{Del} (and the natural relations between heat kernels of
discrete and continuous time Markov chains) to deduce
%
\begin{eqnarray}
c_5 t^{-d/2}\exp\bigl(-c_6
|x-y|^2/t\bigr)&\le& q^n_t(x,y) \le
c_7 t^{-d/2}\exp\bigl(-c_8
|x-y|^2/t\bigr),\nonumber
\\
\eqntext{|x-y| \le t \le\operatorname{diam} ( D),}
\end{eqnarray}
where $q_t^n(x,y)$ is defined as $n^{-d}P^{\omega}_x(X^n_t=y)$ for the
VSRW and $n^{-d}\mu(ny)^{-1}\* P^{\omega}_x(X^n_t=y)$ for the CSRW. Given
these heat kernel estimates, it is then straightforward to verify the
conditions required for the quenched invariance principle by applying
similar arguments to those of Sections~\ref{percsec} and~\ref{boxsec}.

\subsection{Uniform elliptic random divergence form in domains}\label{sec4.3}\label{dfsec}

In this section, we explain how Theorem~\ref{teosilvs} can be applied
in the random divergence form setting. Let $D$ be a uniform domain in
$\mathbb{R}^d$, $d\geq2$.
Assume that we have a random divergence form as follows. There exists a
probability space $(\Omega,{\mathbb P})$ with shift operators $(\tau
_x)_{x\in{\mathbb R}^d}$ that are ergodic, and a symmetric $d\times d$
matrix $A^\omega(x)$ for each $x\in\mathbb{R}^d$ and $\omega\in
\Omega$
such that $A^\omega(x)= A^{\tau_x \omega}(0)$ and
\[
\mathbb P\bigl(c_1I\le A^\omega(x) \le c_2 I
\bigr)=1\qquad\forall x\in\mathbb{R}^d,
\]
where $c_1,c_2\in(0,\infty)$ are deterministic constants. For $n\ge
1$, let
\[
\mathcal{E}^n (f,f)=\frac{1}2\int_{nD}
\nabla f(x) A^\omega(x)\nabla f(x)\,dx.
\]
Let $(Y_t^{n})_{t\ge0}$ be the reflected diffusion process on $n
\overline{D}$
associated with the regular Dirichlet form
$(\mathcal{E}^n, W^{1,2}(nD))$ on $L^2(nD; dx)$,
and set $X^n_t:=n^{-1}Y^{n}_{n^2t}$.
[A natural setting would be to take $D$ to be a cone. In this case $nD=D$,
so the random diffusion matrix $A^\omega(x)$ only needs to be defined
for $x\in D$ rather than for \mbox{$x\in\mathbb{R}^d$}.]
Observe that process $X^n$ takes value in $\overline{D}$.
It is then the case that $X^n=(X_t^n)_{t\geq0}$ converges as
$n\rightarrow\infty$ to a reflecting Brownian motion on $D$ with some
strictly positive covariance matrix $B$, for $\mathbb{P}$-almost-every
realization of the random environment $\omega$.
(Note that $B$ is determined by the invariance principle on the whole
space $\mathbb{R}^d$.)
Indeed, the Dirichlet form of $X^n$ on $L^2(D; dx)$ is
\[
n^{2-d}\mathcal{E}^n(f_n,f_n)=
\frac{1}2\int_D\nabla f(x)A(nx)\nabla f(x)\,dx,
\]
where $f_n(x):=f(n^{-1}x)$.
In view of Section~\ref{sec2.1}, the transition density function of $X^n$ has estimates
(\ref{e23n}) with constants $c_1,\ldots, c_4$ independent of $n$.
In this case, the quenched invariance principle on $\mathbb{R}^d$ is
proved in
\cite{Osa}. (To be precise, in the paper the author assumed $C^2$
smoothness for the coefficients.
However, this was to apply the It\^o formula, and could be avoided by
using the Fukushima decomposition instead.)
Given these, one can easily verify the conditions required for the
quenched invariance principle in $D$. (Note that because of the uniform
ellipticity, condition~(iii) in Theorem~\ref{teosilvs} is trivial in
this case. Moreover, one can extend
the quenched invariance principle of \cite{Osa} to arbitrary starting
points by applying the argument of Theorem~\ref{qip}.) Thus, for
$\mathbb{P}$-almost-every realization of the random environment
$\omega
$ and for every starting point $x\in\overline{D}$, the reflecting
diffusion $X^n$ converges weakly to a reflecting Brownian motion on
$D$. This gives an affirmative answer to the open problem of \cite
{Rho}, pages 1004--1005.

As we mentioned briefly in the \hyperref[sec1]{Introduction}, homogenization of
reflected SDE/PDE on half-planes has been studied for periodic
coefficients in \cite{BLLS,BLP,Tana}, etc., and for random divergence
forms in \cite{Rho}. (Note that their equations contain additional
reflection terms, though the precise framework varies in each paper.)
Homogenization for random divergence forms without extra reflection
terms on bounded $C^2$ domains is discussed in \cite{KLO},
Section~14.4. Although we can only handle symmetric cases,
our results hold for general uniform domains.

\setcounter{equation}{0}
\begin{appendix}
\section*{Appendix}\label{sec5}
\subsection{Proofs for percolation estimates}\label{sec5.1}

The aim of this section is to verify the percolation estimates stated
as Lemmas~\ref{holelem}--\ref{ballcomp},~\ref{meascomp} and~\ref{boxclust}.

For\vspace*{1.5pt} the purpose of proving Lemma~\ref{holelem}, it will be useful to
note that for large $K$ the collection of edges $\widetilde{\mathcal{O}}_2$
stochastically dominates $\widetilde{\mathcal{O}}_3$, the edges of a bond
percolation process on $E_{\mathbb{Z}^d}$ with probability
$p_3=p_3(K)$, where the parameter $p_3$ can be chosen to satisfy $\lim
_{K\rightarrow\infty}p_3=p_1$ (see \cite{ABDH}, Proposition 2.2). In
fact, the proof of this result from \cite{ABDH} further shows that, for
a given value of $K$ (i.e., suitably large), it\vspace*{1pt} is possible to couple
all the relevant random variables in such a way that $\widetilde
{\mathcal
{O}}_3\subseteq\widetilde{\mathcal{O}}_2$ almost-surely. We will
henceforth assume that this is the case, where $K$ is fixed large
enough to ensure that $p_3>p_c^{\mathrm{bond}}(\mathbb{Z}^d)$. We will
also define $\mathcal{O}_3:=\widetilde{\mathcal{O}}_3\cap E_{\mathbb{L}}$
and $\mathcal{C}_3:=\mathcal{C}_\infty(\mathbb{L},\mathcal{O}_3)$. Note
that $\mathcal{C}_3$ is nonempty by the uniqueness of infinite
supercritical bond percolation clusters on $\mathbb{L}$.

\begin{pf*}{Proof of Lemma~\ref{holelem}}
First observe that $\mathcal{O}_3\subseteq\mathcal{O}_2\subseteq
\mathcal{O}_1$. It follows that there exists an infinite connected
component $\mathcal{C}$ of $(\mathbb{L},\mathcal{O}_2)$ such that
$\mathcal{C}_3\subseteq\mathcal{C}\subseteq\mathcal{C}_1$. For
such a
$\mathcal{C}$ (at the moment, we do not know its uniqueness), we have
that $\mathcal{C}_1\setminus\mathcal{C}\subseteq\mathcal
{H}_3\subseteq\mathbb{L}\setminus\mathcal{C}_3$.
Denote by $\mathcal{G}(x)$
the connected component of $\mathbb{L}\setminus\mathcal{C}_3$
containing $x$ [if~$x\in\mathcal{C}_3$, we set $\mathcal
{G}(x)=\varnothing$].
To prove part (i) of the lemma, it suffice to show
that there exist constants $c_1,c_2$ such that: for each $x\in\mathbb{L}$,
%
\begin{equation}
\label{target} \mathbb{P} \bigl(\operatorname{diam} \bigl( \mathcal
{G}(x) \bigr) \geq n
\bigr)\leq c_1e^{-c_2n}.
\end{equation}
For this, we will follow the renormalization argument
used in the proof of
\cite{BiskP}, Proposition 2.3, making the adaptations necessary to deal
with the boundary issues that arise in our setting.

We start by coupling a finite range-dependent site percolation model
with our bond percolation process. For $L\in\mathbb{N}$, $x\in
\mathbb
{Z}^d$, define
\begin{eqnarray*}
Q_L(x)&:=&L(x+e_1+\cdots+e_{d_1})+[0,L]^d
\cap\mathbb{Z}^d,
\\
\widetilde{Q}_{3L}(x) &:=& L(x+e_1+
\cdots+e_{d_1})+[-L,2L]^d\cap\mathbb{Z}^d,
\end{eqnarray*}
where $e_1,\ldots,e_d$ are the standard basis vectors for $\mathbb
{Z}^d$. Given these sets,
let $G_L(x)$
be the event such that:
\begin{itemize}
\item there\vspace*{1pt} exists an $\widetilde{\mathcal{O}}_3$-crossing cluster for
$Q_L(x)$ in $\widetilde{Q}_{3L}(x)$, that is, there is a \mbox{$\widetilde
{\mathcal{O}}_3$-}connected cluster in $\widetilde{Q}_{3L}(x)$ that, for all $d$
directions, joins the `left face' to the `right face' of $Q_L(x)$,\vspace*{1pt}
\item all paths along edges of $\widetilde{\mathcal{O}}_3$ that are
contained in $\widetilde{Q}_{3L}(x)$ and have diameter greater than $L$ are
connected to the (necessarily unique) crossing cluster.
\end{itemize}
We will say that $x\in\mathbb{Z}^d$ is ``white'' if $G_L(x)$ holds and
``black'' otherwise, and make the important observation that if two
neighboring vertices are white, then their crossing clusters must be
connected. Since $p_3>p_c^{\mathrm{bond}}(\mathbb{Z}^d)$, we can apply
(the bond percolation version of) \cite{PP}, Theorem~5, for $d=2$ and
\cite{Pisz}, Theorem~3.1, for $d\geq3$ to deduce that
%
\begin{equation}
\label{glim} \lim_{L\rightarrow\infty}\mathbb{P} \bigl(G_L(x)
\bigr)=1
\end{equation}
(cf. \cite{Ant-P}, (2.24)). Moreover, although $(\mathbf
{1}_{G_L(x)})_{x\in\mathbb{Z}^d}$ are not
independent,
the dependence between these random variables
is of finite range. Thus, by \cite{LSS}, Theorem~0.0,
one can suppose that, for suitably large $L$, the collection $(\mathbf
{1}_{G_L(x)})_{x\in\mathbb{Z}^d}$ dominates a site percolation process
on $\mathbb{Z}^d$ of density arbitrarily close to 1. Noting that for
any infinite connected graph $G$ with maximal vertex degree $\Delta$
the critical site and bond percolation probabilities satisfy
\[
p_c^{\mathrm{site}}(G)\leq1-\bigl(1-p_c^{\mathrm{bond}}(G)
\bigr)^{\Delta-1}
\]
(\cite{GS}, Theorem~3), we have that $p_c^{\mathrm{site}}(\mathbb{L})$
is bounded above by $1-(1-p_c^{\mathrm{bond}}(\mathbb{Z}^d))^{2d-1}<1$.
Hence, by taking $L$ suitably large, it is possible to assume that
there is a nonzero probability of 0 being contained in an infinite
cluster of white vertices in $\mathbb{L}$. From this result, a standard
ergodicity argument with respect to the shift $x\mapsto x+e_1+\cdots+e_{d_1}$
allows one to check that,
$\mathbb{P}$-a.s., there exists at least one infinite connected cluster
of white sites in $\mathbb{L}$. In particular, writing $\mathcal{D}(x)$
for the connected cluster of white sites containing a particular vertex
$x\in\mathbb{L}$ [taking $\mathcal{D}(x):=\varnothing$ if $G_L(x)$ does
not occur], we obtain that the set
%
\begin{equation}
\label{dinfdef} \mathcal{D}_\infty:= \bigl\{x\in\mathbb{L}\dvtx\mathcal{D}(x)\mbox{ is infinite} \bigr\}
\end{equation}
is nonempty, $\mathbb{P}$-a.s.

Let $\mathcal{C}(x)$ be the connected component of $\mathbb
{L}\setminus\mathcal{D}_\infty$ containing $x$ (we set this to be the
empty set if $x\in\mathcal{D}_\infty$).
The next step of the proof is to check that: for $x\in\mathbb{L}$,
%
\begin{equation}
\label{cbound} \mathbb{P} \bigl(\operatorname{diam} \bigl(\mathcal
{C}(x)\bigr) \geq n
\bigr)\leq c_3e^{-c_4n}.
\end{equation}
To do this, we start by introducing some notions of set boundaries that
will be useful. For $x\notin\mathcal{D}_\infty$, the inner boundary of
$\mathcal{C}(x)$ is the set
\[
\partial^{\mathrm{in}}\mathcal{C}(x):= \bigl\{y\in\mathcal{C}(x)\dvtx
y\mbox{ is adjacent to a vertex in $\mathbb{L}\setminus\mathcal{C}(x)$} \bigr\}.
\]
It is simple to check from its construction that all the vertices in
this set are black. Since $\mathbb{L}\setminus\mathcal
{C}(x)\supseteq
\mathcal{D}_\infty\neq\varnothing$, then $\mathbb{L}\setminus
\mathcal
{C}(x)$ contains at least one infinite connected component, $\mathcal
{D}$ say. The outer boundary of $\mathcal{D}$ is given by
%
\begin{equation}
\label{outer} \partial^{\mathrm{out}}\mathcal{D}:= \{y\in\mathbb
{L}\setminus
\mathcal{D}\dvtx y\mbox{ is adjacent to a vertex in $\mathcal{D}$} \}.
\end{equation}
With $\mathcal{D}$ being disjoint from $\mathcal{C}(x)$, we can also
define the part of its outer boundary visible from $\mathcal{C}(x)$ by setting
%
\begin{eqnarray}\label{vis}
\partial^{\mathrm{vis}\mathcal{C}(x)}\mathcal{D} &:=& \bigl\{
y\in\partial ^{\mathrm{out}}\mathcal{D}: \mbox{there exists a path from $y$ to $
\mathcal{C}(x)$}
\nonumber\\[-8pt]\\[-8pt]
&&\hspace*{83pt}\mbox{in $\mathbb{L}$ that is disjoint from $\mathcal{D}$}
\bigr\}.\nonumber
\end{eqnarray}

We claim the following relationship between the various boundary sets:
%
\begin{equation}
\label{incclaim} \partial^{\mathrm{vis}\mathcal{C}(x)}\mathcal
{D}=\partial^{\mathrm
{out}}
\mathcal{D}\subseteq\partial^{\mathrm{in}}\mathcal{C}(x).
\end{equation}
To verify the equality, first suppose that there exists a vertex $y\in
\partial^{\mathrm{out}}\mathcal{D}\setminus\mathcal{C}(x)$, then
$y\in
\mathbb{L}\setminus\mathcal{C}(x)$ and we can find a vertex $z\in
\mathcal{D}\subseteq\mathbb{L}\setminus\mathcal{C}(x)$ such that $y$
and $z$ are adjacent. This implies that $y$ and $z$ are in the same
connected component of $\mathbb{L}\setminus\mathcal{C}(x)$, which
is a
contradiction because $y\notin\mathcal{D}$ by definition. Hence
$\partial^{\mathrm{out}}\mathcal{D}\subseteq\mathcal{C}(x)$, and so
[noting that $\mathcal{C}(x)\setminus\mathcal{D}=\mathcal{C}(x)$]
\[
\partial^{\mathrm{out}}\mathcal{D}= \bigl\{y\in\mathcal{C}(x)\dvtx
y\mbox{ is
adjacent to a vertex in $\mathcal{D}$} \bigr\}\subseteq\partial
^{\mathrm{vis}\mathcal{C}(x)}
\mathcal{D}.
\]
Since the opposite inclusion is trivial, we obtain the equality at
(\ref
{incclaim}). From $\partial^{\mathrm{out}}\mathcal{D}\subseteq
\mathcal
{C}(x)$, the inclusion at (\ref{incclaim}) is also clear.

We proceed by applying the conclusion of the previous paragraph to show
that $\mathcal{D}=\mathbb{L}\setminus\mathcal{C}(x)$. First, the
boundary connectivity result of \cite{Timar}, Lemma 2, implies that
$\partial^{\mathrm{vis}\mathcal{C}(x)}\mathcal{D}$ is $*$-connected.
Combining this with (\ref{incclaim}), we obtain that $\partial
^{\mathrm
{out}}\mathcal{D}$ is a $*$-connected set of black vertices [recall
that the vertices of $\partial^{\mathrm{in}}\mathcal{C}(x)$ are black].
Secondly, note that if $\mathbb{P}_p$ is the law of a parameter $p$
site percolation process on $\mathbb{Z}_d$ and $\mathcal{C}^*$ is the
corresponding $*$-connected component of closed vertices containing~$0$, then for suitably large $p$ we have that
%
\begin{equation}
\label{exptails} \mathbb{P}_p \bigl(\mathcal{C}^*\geq n \bigr)\leq
c_5e^{-c_6n}
\end{equation}
(see \cite{AiB}, Theorem 7.3 and \cite{AiN}, Proposition 7.6). In
particular, it is easy to check from this that all $*$-connected
components of closed vertices in the site percolation process with this
choice of $p$ are finite, $\mathbb{P}_p$-a.s. Hence, because $(\mathbf
{1}_{G_L(x)})_{x\in\mathbb{Z}^d}$ dominates a site percolation process
whose parameter can be made arbitrarily close to~1 by taking $L$
suitably large, it must be the case that, for large $L$, $\partial
^{\mathrm{out}}\mathcal{D}$ is $\mathbb{P}$-a.s. a finite set. Since
$\mathcal{D}$ is infinite, it readily follows that $\mathbb
{L}\setminus
\mathcal{D}$ is also finite. Now, suppose $\mathcal{D}_1$ is a
connected component of $\mathbb{L}\setminus\mathcal{C}(x)$ distinct
from $\mathcal{D}$ and such that $\mathcal{D}_1\cap\mathcal
{D}_\infty
\neq\varnothing$. By the definition of $\mathcal{D}_\infty$, it holds
that $\mathcal{D}_1\cap\mathcal{D}(y)\neq\varnothing$ for some $y$ such
that $\mathcal{D}(y)$ is an infinite set. Since $\mathcal{D}_1\cup
\mathcal{D}(y)$ is an infinite connected component of $\mathbb
{L}\setminus\mathcal{C}(x)$, it must be the case that $\mathcal{D}_1$
is infinite. However, this contradicts the finiteness of $\mathbb
{L}\setminus\mathcal{D}$, and so no such $\mathcal{D}_1$ can exist.
Thus, it must be the case that if $\mathcal{D}_1$ is a connected
component of $\mathbb{L}\setminus\mathcal{C}(x)$ distinct from
$\mathcal{D}$, then $\mathcal{D}_1\cap\mathcal{D}_\infty=\varnothing$.
Since, similar to the results of the previous paragraph, we have that
$\partial^{\mathrm{out}}\mathcal{D}_1\subseteq\partial^{\mathrm
{in}}\mathcal{C}(x)$, the set $\mathcal{D}_1\cup\mathcal{C}(x)$
must be
a connected component of $\mathbb{L} \setminus\mathcal{D}_\infty$. By
the definition of $\mathcal{C}(x)$, this implies that $\mathcal
{D}_1=\varnothing$, which yields
$\mathcal{D}=\mathbb{L}\setminus\mathcal{C}(x)$ as required.

An immediate corollary of the equality $\mathcal{D}=\mathbb
{L}\setminus
\mathcal{C}(x)$ is that $\partial^{\mathrm{in}}\mathcal
{C}(x)=\partial
^{\mathrm{out}}\mathcal{D}$,
which as we have already established, is a finite $*$-connected
component of black vertices. We will use these results to finally prove
(\ref{cbound}). Note first that $\operatorname{diam}(\mathcal
{C}(x))=\operatorname{diam} (\partial^{\mathrm{in}}\mathcal{C}(x))$.
Hence, writing
$\partial B(x,m)$ for the vertices of $\mathbb{L}$ at an $\ell_\infty$
distance $m$ from $x$,
\begin{eqnarray*}
\mathbb{P} \bigl(\operatorname{diam} \bigl(\mathcal{C}(x)\bigr)\geq n
\bigr)&\leq&
\sum_{m=0}^\infty\mathbb{P} \bigl(
\operatorname{diam} \bigl(\partial^{\mathrm
{in}}\mathcal{C}(x)\bigr) \geq n,
\partial^{\mathrm{in}}\mathcal{C}(x)\cap\partial B(x,m)\neq
\varnothing\bigr)
\\
&\leq& \sum_{m=0}^\infty\sum
_{y\in\partial B(x,m)}\mathbb{P} \bigl(\operatorname{diam} \bigl
(\mathcal{C}^*(y)
\bigr) \geq n\vee m \bigr),
\end{eqnarray*}
where $\mathcal{C}^*(y)$ is the $*$-connected component of black
vertices in $\mathbb{Z}^d$ containing $y$. By again comparing
$(\mathbf
{1}_{G_L(x)})_{x\in\mathbb{Z}^d}$ to a site percolation process, it is
possible to apply (\ref{exptails}) to deduce that the tail of the
probability in the above sum is bounded above by $c_5e^{-c_6(n\vee
m)}$. The estimate at (\ref{cbound}) follows.

We now return to the problem of deriving the estimate at (\ref
{target}). For $x\in\mathbb{Z}^d$, define the set
$Q_L'(x):=L(x+e_1+\cdots+e_{d_1})+[0,L)^d\cap\mathbb{Z}^d$, so that
$(Q_L'(x))_{x\in\mathbb{Z}^d}$ is a partition of $\mathbb{{Z}}^d$. For
$x\in\mathbb{L}$, let $a(x)$ be the closest element of $\mathbb{L}$,
with respect to $\ell_1$ distance, to the $x'\in\mathbb{Z}^d$ such
that $x\in Q_L'(x')$. [Only when $x$ is within a distance $L$ of the
boundary of $\mathbb{L}$ does $x'\neq a(x)$.] It is then easy to check
that if $\operatorname{diam} (\mathcal{G}(x)) \geq L$, then
%
\begin{equation}
\label{gcontain} \mathcal{G}(x)\subseteq\bigcup_{x'\in\mathcal{C}(a(x))}
\widetilde{Q}_{3L}\bigl(x'\bigr)
\end{equation}
(cf. \cite{BiskP}, (3.7)). In particular, this implies that, if
$\operatorname{diam} (\mathcal{G}(x)) \geq L$, then it must be the case that
$\operatorname{diam} (\mathcal{G}(x)) \leq3L \operatorname{diam}
(\mathcal{C}(a(x)))$.
Consequently, (\ref{target}) follows from~(\ref{cbound}), and thus the
proof of part (i) is complete.

Given part (i), we observe
\begin{eqnarray*}
&& \mathbb{P} \bigl(\exists x\in[-n,n]^d\cap\mathbb{L}\dvtx \#
\mathcal{H}(x)\geq(\log n)^{d+1} \bigr)
\\
&&\qquad \leq (2n+1)^d\sup_{x\in\mathbb{L}}\mathbb{P} \bigl( x\in
\mathcal{C}_1, \#\mathcal{H}(x)\geq(\log n)^{d+1} \bigr)
\\
&&\qquad \leq(2n+1)^d\sup_{x\in\mathbb{L}}\mathbb{P} \bigl( x\in
\mathcal{C}_1, \operatorname{diam} \bigl(\mathcal{H}(x)\bigr)\geq(\log
n)^{(d+1)/d} \bigr)
\\
&&\qquad \leq(2n+1)^dc_1e^{-c_2(\log n)^{(d+1)/d}}.
\end{eqnarray*}
Since this is summable in $n$, part (ii) follows by a Borel--Cantelli argument.
\end{pf*}

We now work toward the proof of Lemma~\ref{distlem}.

\begin{pf*}{Proof of Lemma~\ref{distlem}}
To establish the
bound in (\ref{upperd1}), let us start by recalling/adapting some
definitions from the previous proof. In particular, for $x\in\mathbb
{Z}^d$, define $Q_L(x)$ and $\widetilde{Q}_{3L}(x)$ as in the proof of
Lemma~\ref{holelem}. Moreover, let $G_L(x)$ be defined similarly, but
with $\widetilde{\mathcal{O}}_3$ replaced by $\widetilde{\mathcal
{O}}_1$, and
redefine $x$ being ``white'' to mean that this version of $G_L(x)$
holds (and say $x$ is ``black'' otherwise). Note that the statement
(\ref{glim}) remains true with this definition of $G_L(x)$, and the
dependence between the random variables $(\mathbf{1}_{G_L(x)})_{x\in
\mathbb{Z}^d}$ is only finite range, and so we can suppose that it
dominates a dominates a site percolation process on $\mathbb{Z}^d$ of
density arbitrarily close to 1.

Now, fix $x,y\in\mathbb{L}$, and recall the definition of $a(x)$ from
the proof of Lemma~\ref{holelem}. If $n$ is the $\ell_1$ distance
between $a(x)$ and $a(y)$, then there exists a nearest neighbor path
$a_0,\ldots,a_n$ in $\mathbb{L}$ such that $a_0=a(x)$ and $a_n=a(y)$. We
claim that if $x$ and $y$ are both contained in $\mathcal{C}_1$, then
there exists a path from $x$ to $y$ along edges of $\mathcal{O}_1$
whose vertices all lie in
%
\begin{equation}
\label{set} \bigcup_{i=0}^n\bigcup
_{b\in\overline{\mathcal{C}}^*({a_i})}\widetilde{Q}_{3L}(b),
\end{equation}
where $\overline{\mathcal{C}}^*(a):=\{a\}$ if $a$ is white,
otherwise $\overline{\mathcal{C}}^*(a):=\mathcal{C}^*(a)\cup\partial
^{\mathrm{out}}\mathcal
{C}^*(a)$, where $\mathcal{C}^*(a)$ is the $*$-connected component of
black sites in $\mathbb{L}$ containing $a$ [$\partial^{\mathrm
{out}}\mathcal{C}^*(a)$ is the outer boundary of $\mathcal{C}^*(a)$,
defined similarly to (\ref{outer})]. This is essentially \cite{Ant-P},
Proposition 3.1, rewritten for $\mathbb{L}$ instead of
$\mathbb
{Z}^d$. The one slight issue with modifying the proof of this result to
our situation is that, unlike the $\mathbb{Z}^d$ case, the outer
boundary in $\mathbb{L}$ of a finite connected cluster of vertices,
$\mathcal{C}$ say, is no longer $*$-connected in general and so it is
not possible to run around it in quite the same way. However, this
problem is readily overcome by applying \cite{Timar}, Lemma 2, which
implies that for each $x\notin\mathcal{C}$, the part of the outer
boundary of $\mathcal{C}$ that is visible from $x$, $\partial
^{\mathrm
{vis}(x)}\mathcal{C}$ [cf. (\ref{vis})], is $*$-connected. A simple
estimate of the number of vertices in the set at (\ref{set}) yields
\[
d_1(x,y)\leq(3L+1)^d\sum_{i=0}^n
\#\overline{\mathcal{C}}^*(a_i)\leq(3L+1)^d\sum
_{i=0}^n \bigl(1+3^d\#{
\mathcal{C}}^*(a_i) \bigr).
\]
[We take ${\mathcal{C}}^*(a):=\varnothing$ if $a$ is white.] Clearly,
$\mathcal{C}^*(a)\subseteq\widetilde{\mathcal{C}}^*(a)$, where
$\widetilde
{\mathcal{C}}^*(a)$ is the $*$-connected component of black sites in
$\mathbb{Z}^d$ containing $a$. Consequently, we obtain that
\[
\mathbb{P} \bigl(x,y\in\mathcal{C}_1\mbox{ and }d_1(x,y)
\geq cR \bigr)\leq\mathbb{P} \Biggl((3L+1)^d\sum
_{i=0}^n \bigl(1+3^d\#\widetilde{{\mathcal
{C}}}^*(a_i) \bigr)\geq cR \Biggr).
\]
Moreover, applying \cite{DeuP}, Lemma 2.3, as in the proof of \cite
{Ant-P}, Theorem 1.1, we may replace the summands on the right-hand
side by independent ones, each with the same distribution as the term
they are replacing. Recalling from (\ref{exptails}) the exponential
bound for the size of a $*$-connected vacant cluster in a site
percolation process of parameter $p$ close to 1, one readily obtains
from this the bound at (\ref{upperd1}). (It is useful to also note that
$n\leq c |x-y|/L$.)

We proceed next with the proof of the second bound. In this direction,
let us begin by defining a metric $d_Z$ on $\mathcal{C}_2$ related to
the process $Z$ introduced in Section~\ref{hkestsec}, that is, the time
change of $Y$ with time in $\mathcal{H}$ cut out. Assume that $K$ is
large enough so that the conclusions of Lemma~\ref{holelem} hold.
Define a set of edges $E_Z'$ by supposing, for $x,y\in\mathcal{C}_2$,
$\{x,y\}\in E_Z'$ if and only if $\{x,y\}\notin\mathcal{O}_2$ and also
there exists a path $x=z_0,z_1,\ldots,z_k=y$ such that $z_1,\ldots,z_{k-1}\in\mathcal{H}$ and $\{z_{i-1},z_i\}\in\mathcal{O}_1$ for
$i=1,\ldots,k$. Thus, the jumps of $Z$ will be on edges in either
$\mathcal{O}_2$ or $E_Z'$. Set $E_Z:=\mathcal{O}_2\cup E_Z'$, and let
$d_Z$ be the graph distance on $(\mathcal{C}_2,E_Z)$. Our first goal
will be to prove that there exist constants $c_1,c_2,c_3$ such that:
for every $x,y\in\mathbb{L}$,
%
\begin{equation}
\label{lower} \mathbb{P} \bigl(x,y\in\mathcal{C}_2\mbox{ and
}d_Z(x,y)\leq c_1^{-1}|x-y| \bigr)\leq
c_2e^{-c_3|x-y|},
\end{equation}
where $|x-y|$ is the Euclidean distance between $x$ and $y$.

For proving (\ref{lower}), we suppose that the definition of $G_L(x)$
reverts to that given in the proof of Lemma~\ref{holelem}, that is, in
terms of $\widetilde{\mathcal{O}}_3$. Also, define $G_L'(x)$ to be the
event that there are no edges of the set $\widetilde{\mathcal
{O}}_1\setminus\widetilde{\mathcal{O}}_3$ connecting two vertices of
$\widetilde{Q}_{3L}(x)$, so that if $G_L(x)\cap G_L'(x)$ holds,\vspace*{1pt} then so do
the defining properties of $G_L(x)$ when~$\widetilde{\mathcal{O}}_3$ is
replaced by $\widetilde{\mathcal{O}}_2$. Clearly, for fixed $L$,
$\mathbb
{P} (G_L'(x)^c )\rightarrow0$ as $p_3\rightarrow p_1$
(i.e., $K\rightarrow\infty$). Hence, for any $\delta$, by first
choosing $L$ and then $K$ large, we can ensure $\mathbb{P}
(G_L(x)\cap G_L'(x) )\geq1-\delta$. For the remainder of this
proof, we redefine $x$ being ``white'' to mean that $G_L(x)\cap
G_L'(x)$ holds, and say $x$ is ``black'' otherwise.

Similarly to above, the finite range dependence of the random variables
in question means that it is possible to suppose that $(\mathbf
{1}_{G_L(x)\cap G_L'(x)})_{x\in\mathbb{Z}^d}$\vspace*{1pt} stochastically dominates
a collection $(\eta(x))_{x\in\mathbb{Z}^d}$ of independent and
identical Bernoulli random variables whose parameter $p$ is arbitrarily
close to 1. Let $\mathcal{C}_\infty$ be the vertices of~$\mathbb{L}$
that are contained in an infinite connected component of $\{x\in
\mathbb
{L}\dvtx\eta(x)=1\}$ [cf.~(\ref{dinfdef})]. By arguments from the
proof of
Lemma~\ref{holelem}, we have that if $p$ is large enough, then this set
is nonempty and its complement in $\mathbb{L}$ consists of finite
connected components, $\mathbb{P}$-a.s. Now, as in the proof of \cite
{BiskP}, Lemma 3.1, we ``wire'' the holes of $\mathcal{C}_\infty$ by
adding edges between every pair of sites that are contained in a
connected component of $\mathbb{L}\setminus\mathcal{C}_\infty$ or its
outer boundary, and denote the induced graph distance by $d'$. By
proceeding almost exactly as in \cite{BiskP}, it is then possible to
show that, for suitably large $L$ and $K$: for $x,y\in\mathbb{L}$,
%
\begin{equation}
\label{dest} \mathbb{P} \bigl(d'(x,y)\leq\tfrac{1}{2}|x-y|
\bigr)\leq e^{-|x-y|}.
\end{equation}
[The one modification needed depends on the observation that, similar
to what was deduced in the proof of Lemma~\ref{holelem}, the inner
boundary of any connected component of $\mathbb{L}\setminus\mathcal
{C}_\infty$ is $*$-connected and consists solely of vertices with
$\eta
(x)=0$.] Finally, a minor adaptation of (\ref{gcontain}) yields, for
$x\in\mathcal{C}_1$ with $\operatorname{diam}(\mathcal{H}(x))\geq L$,
\[
\mathcal{H}(x)\subseteq\bigcup_{x'\in\mathcal{C}(a(x))}
\widetilde{Q}_{3L}\bigl(x'\bigr),
\]
where $\mathcal{C}(a(x))$ is now the connected component of $\mathbb
{L}\setminus\mathcal{C}_\infty$ containing $a(x)$. It follows that if
$x,y\in\mathcal{C}_2$, then $d_Z(x,y)\geq d'(a(x),a(y))$, cf. \cite
{BiskP}, (3.10). Therefore, since it also holds for $|x-y|\geq3L$ that
$|a(x)-a(y)|\geq c|x-y|/L$, the bound at~(\ref{lower}) can be obtained
from (\ref{dest}).

Finally, note that, since $\mu_e\in[K^{-1},K]$ for every $e\in
\mathcal
{O}_1$ such that $e\cap e'\neq\varnothing$ for some $e'\in\mathcal
{O}_2$, it holds that $t(e)\in[C_A\wedge K^{-1/2}, K^{1/2}]$ for such
edges. Moreover, for every other $e\in\mathcal{O}_1$, we have
$t(e)\geq0$. As a consequence, the metric $\bar{d}_1$ is
bounded below
by a constant multiple of ${d}_Z$ on $\mathcal{C}_2$. Applying this, as
well as setting $\partial^{\mathrm{out}}\mathcal{H}(x)=\{x\}$ for
$x\notin
\mathcal{H}$, it follows that
\begin{eqnarray*}
&& \mathbb{P} \bigl(x,y\in\mathcal{C}_1\mbox{ and }\bar
{d}_1(x,y)\leq c_4^{-1}|x-y| \bigr)
\\
&&\qquad \leq\mathbb{P} \Bigl(x,y\in\mathcal{C}_1\mbox{ and }\inf
_{x'\in
\partial^{\mathrm{out}}\mathcal{H}(x), y'\in\partial^{\mathrm
{out}}\mathcal
{H}(y)}{d}_Z\bigl(x',y'
\bigr)\leq c_5^{-1}|x-y| \Bigr)
\\
&&\qquad \leq\sum_{x',y'\dvtx |x-x'|, |y-y'|\leq|x-y|/4} \mathbb{P} \bigl(x',y'
\in\mathcal{C}_2\mbox{ and } {d}_Z\bigl(x',y'
\bigr)\leq c_5^{-1}|x-y| \bigr)
\\
&&\quad\qquad{}+2\mathbb{P} \bigl(x\in\mathcal{C}_1\mbox{ and } \operatorname
{diam} \bigl(\mathcal{H}(x)\bigr)+ 1\geq|x-y|/4 \bigr).
\end{eqnarray*}
From this, the bound at (\ref{lowerd1}) is a straightforward
consequence of Lemma~\ref{holelem}(i) and (\ref{lower}).
\end{pf*}

Given Lemma~\ref{distlem}, the proof of Lemma~\ref{ballcomp} is
straightforward.

\begin{pf*}{Proof of Lemma~\ref{ballcomp}}
Since $\bar{d}_1\leq C_Ad_1$, the inclusion $B_1(x,R)\subseteq\break 
\overline{B}_1(x,C_A R)$ always holds. We will thus concern ourselves
with the
other two inclusions only. First, by the inequality at (\ref{lowerd1}),
\begin{eqnarray*}
&& \mathbb{P} \bigl(x\in\mathcal{C}_1\mbox{ and }\overline
{B}_1(x,C_AR)\nsubseteq B_E(x,c_2R)
\bigr)
\\
&&\qquad \leq \sum_{y\notin B_E(x,c_2R)}\mathbb{P} \bigl(x,y\in\mathcal
{C}_1\mbox{ and }\bar{d}_1(x,y)\leq
C_AR \bigr)
\\
&&\qquad \leq c_5e^{-c_6R}
\end{eqnarray*}
for suitably large $c_2$. Second,
\begin{eqnarray*}
&& \mathbb{P} \bigl(x\in\mathcal{C}_1\mbox{ and }\mathcal{C}_1
\cap{B}_E(x,c_1R)\nsubseteq B_1(x,R) \bigr)
\\
&&\qquad
\leq\sum_{y\in B_E(x,c_1R)}\mathbb{P} \bigl(x,y\in\mathcal
{C}_1\mbox{ and } {d}_1(x,y)\geq R \bigr)
\end{eqnarray*}
and applying (\ref{upperd1}) yields a bound of the form $c_7e^{-c_8R}$,
thereby completing the proof.
\end{pf*}

Next, the comparison of measures stated as Lemma~\ref{meascomp}.

\begin{pf*}{Proof of Lemma~\ref{meascomp}}
First note that any point $x\in Q$
that is contained in $\widetilde{\mathcal{C}}_1\setminus\mathcal{C}_1$
must lie in a connected component of $(Q\setminus\mathcal
{C}_1,\mathcal{O}_1)$ that meets the inner boundary of $Q$, which we
denote here by $\partial^{\mathrm{in}}Q$. Moreover, we recall that any
points $x\in Q$ that are contained in $\mathcal{C}_1$ must also be
contained in $\widetilde{\mathcal{C}}_1$. It follows that
\[
\tilde{\nu}(Q)-{\nu}(Q)\leq\sum_{x\in\partial^{\mathrm{in}}Q}\#
\mathcal{F}(x),
\]
where $\mathcal{F}(x)$ is the connected component of $\mathbb
{L}\setminus\mathcal{C}_1$ containing $x$. Now, similar to (\ref
{target}), we have that
\[
\mathbb{P} \bigl(\operatorname{diam} \bigl(\mathcal{F}(x)\bigr)\geq n
\bigr)\leq
c_1e^{-c_2n},
\]
uniformly in $x\in\mathbb{L}$. Since $\#\partial^{\mathrm{in}}Q$ is bounded
above by $c_3n^{d-1}$, the lemma follows.
\end{pf*}

Finally, we prove Lemma~\ref{boxclust}.

\begin{pf*}{Proof of Lemma~\ref{boxclust}}
First, observe that
\begin{eqnarray*}
&& \mathbb{P} \bigl(\mathcal{C}^+(Q_2)\nsubseteq\mathcal
{C}_1 \bigr)
\\
&&\qquad \leq\mathbb{P} \bigl(\mathcal{C}^+(Q_2)\leq\varepsilon
|Q_2| \bigr)+\mathbb{P} \bigl(\operatorname{diam} \bigl(\mathcal{C}^+(Q_2)
\bigr) \geq\varepsilon^{1/d}(2n+1), \mathcal{C}^+(Q_2)
\nsubseteq\mathcal{C}_1 \bigr).
\end{eqnarray*}
As in the proof of Proposition~\ref{teoVNG-lem}, the first term here
is bounded above by $c_1e^{-c_2n}$. The second term is bounded above by
\begin{eqnarray*}
&& c_3 n^d \sup_{x\in Q_2}\mathbb{P} \bigl(
\mbox{the diameter of the connected component}
\\[-3pt]
&&\hspace*{93pt}\mbox{of $\mathbb{Z}_+^d
\setminus\mathcal{C}_1$ containing $x$ is }\geq c_4n
\bigr).
\end{eqnarray*}
That this admits a bound of the form $c_5e^{-c_6n}$ follows from (\ref
{target}) (replacing $\mathcal{C}_3$ by~$\mathcal{C}_1$).

Consequently, to complete the proof, it will suffice to show that
\[
\mathbb{P} \bigl(\mathcal{C}^+(Q_2)\cap Q_1\subset
\mathcal{C}_1\cap Q_1 \bigr)\leq c_7e^{-c_8n}.
\]
For the event on the left-hand side of the above to hold, it must be
the case that there exists an open path in $(Q_2,\mathcal{O}_1)$ of
diameter at least $(1-\varepsilon)n$ that is not part of $\mathcal
{C}^+({Q}_2)$. Moreover, as in the previous paragraph, we have that,
with probability at least $1-c_1e^{-c_2n}$, $\operatorname{diam}
(\mathcal
{C}^+(Q_2)) \geq\varepsilon^{1/d}(2n+1)$. However, by (the bond
percolation version of) \cite{PP}, Theorem 5, we have that with
probability at least $1-c_9e^{-c_{10}n}$, there is a unique open
cluster in $(Q_2,\mathcal{O}_1)$ of diameter at least $\varepsilon
^{1/d}(2n+1)$. Hence, by taking $\varepsilon$ suitably small, the
result follows.
\end{pf*}

\subsection{Arbitrary starting point quenched invariance principle}\label{sec5.2}

This section contains the proof of Theorem~\ref{qip}. For it, we note
that the full $\mathbb{Z}^d$ model also satisfies the conclusions of
Propositions~\ref{teotightABDH} and~\ref{teobiskpreLem}, as well as
Assumption~\ref{teoassumpHC} (in the same sense as we checked for the
$\mathbb{L}$ model in Section~\ref{qipsec}).

\begin{pf*}{Proof of Theorem~\ref{qip}}
To begin with, we recall the quenched invariance principle of \cite
{ABDH}, Theorem 1.1(a), for the VSRW started at the origin: there exists
a deterministic constant $c>0$ such that, for $\mathbb
{P}_1$-a.e. $\omega
$, the laws of the processes $\widetilde{Y}^n$ under $\widetilde
{P}_{0}^\omega$
converge weakly to the law of $(B_{ct})_{t\geq0}$, where $(B_t)_{t\geq
0}$ is standard Brownian motion on $\mathbb{R}^d$ started from $0$.
Here, $\mathbb{P}_1$ is the conditional law $\mathbb{P}(\cdot| 0\in
\widetilde{\mathcal{C}}_1)$. Moreover, we note that by proceeding as in
\cite{ABDH}, Remark 5.16, one can check that the result remains true if
$\mathbb{P}_1$ is replaced by $\mathbb{P}$, and $\widetilde
{P}_{0}^\omega$
is replaced by $\widetilde{P}_{x_0}^\omega$, where $x_0$ is chosen to be
the (not\vspace*{1pt} necessarily uniquely defined) closest point to the origin in
the infinite cluster $\widetilde{\mathcal{C}}_1$.

Given the above result, the argument of this paragraph applies for
$\mathbb{P}$-a.e. $\omega$. Fix $x\in\mathbb{R}^d$, $\varepsilon>0$,
and define
\begin{eqnarray*}
\tau^n &:=&\inf\bigl\{t>0\dvtx \widetilde{Y}^n_t\in
B_E(x,\varepsilon) \bigr\},
\\
\tau^B&:=&\inf\bigl\{t>0\dvtx B_{ct}\in B_E(x,
\varepsilon) \bigr\}.
\end{eqnarray*}
A standard argument gives us that, jointly with the convergence of the
previous paragraph, $\tau^n$ converges in distribution to $\tau^B$.
Letting $T>0$ be a deterministic constant, it follows that the laws of
the processes $(\widetilde{Y}^n_{\tau^n+t})_{t\geq0}$ under $\widetilde
{P}_{x_0}^\omega(\cdot| \tau^n\leq T)$ converge weakly to the law of
$(B_{c(\tau+t)})_{t\geq0}$, started from $0$ and conditional on \mbox{$\tau
\leq T$}. Consequently, the Markov property gives us that for every
bounded, continuous function $f\dvtx C([0,\infty),\mathbb
{R}^d)\rightarrow
\mathbb{R}$,
%
\begin{eqnarray}
\label{c1} && \int_{\overline{B}_E(x,\varepsilon)}\widetilde{E}^\omega_{ny}
f\bigl(\widetilde{Y}^n\bigr) P^\omega_{x_0} \bigl(
\widetilde{Y}^n_{\tau^n}\in dy|\tau^n\leq T \bigr)
\nonumber\\[-8pt]\\[-8pt]
&&\qquad \rightarrow\int_{\overline{B}_E(x,\varepsilon)}{E}^B_y
f(B_{c\cdot}) P_{0}^B \bigl(B_{c\tau^B}\in
dy|\tau^B\leq T \bigr),\nonumber
\end{eqnarray}
where $\overline{B}_E(x,\varepsilon)$ is the closure of
$B_E(x,\varepsilon
)$, $P^B_x$ is the law of the standard Brownian motion $B$ started from
$x$, and $E^B_x$ is the corresponding expectation. Furthermore, it is
elementary to check for such $f$ that, as $\varepsilon\rightarrow0$,
%
\begin{equation}
\label{c2} \int_{\overline{B}_E(x,\varepsilon)}{E}^B_y
f(B_{c\cdot}) P_{0}^B \bigl(B_{c\tau^B}\in
dy|\tau^B\leq T \bigr)\rightarrow{E}^B_x
f(B_{c\cdot}).
\end{equation}
Suppose\vspace*{2pt} that the following also holds for every sequence of starting
points $x_n\in n^{-1}\widetilde{\mathcal{C}}_1$ such that
$x_n\rightarrow
x$, every finite collection of times $0<t_1<t_2<\cdots<t_k$ and all
bounded, continuous functions $f_i\dvtx\mathbb{R}^d\rightarrow\mathbb{R}$,
$i=1,\ldots,k$,
%
\begin{equation}
\label{c3} \lim_{\varepsilon\rightarrow0}\limsup_{n\rightarrow\infty}
\biggl\llvert\int_{\overline{B}_E(x,\varepsilon)}\widetilde{E}^\omega_{ny}
f\bigl(\widetilde{Y}^n_{\mathbf
{t}}\bigr) P^\omega_{x_0}
\bigl(\widetilde Y^n_{\tau^n}\in dy|\tau^n\leq T
\bigr) -\widetilde{E}^\omega_{nx_n} f\bigl(\widetilde{Y}^n_{\mathbf{t}}
\bigr)\biggr\rrvert=0,\hspace*{-40pt}
\end{equation}
where $\widetilde{Y}^n_{\mathbf{t}}:=(\widetilde{Y}^n_{{t_1}},\widetilde
{Y}^n_{t_2},\ldots,\widetilde{Y}^n_{t_k})$. Combining (\ref{c1}), (\ref
{c2}) and (\ref{c3})
readily yields that the finite-dimensional distributions of $\widetilde
{Y}^n$ (under $P^\omega_{nx_n}$) converge to those of $B_{c\cdot}$
(under $P^B_x$).
Moreover, from the full plane version of Proposition~\ref{teotightABDH}, we have that the laws of $\widetilde{Y}^n$ under
$P^\omega
_{nx_n}$ are tight. These two facts yields the desired quenched
invariance principle.

To complete the argument of the previous paragraph and the proof of the
theorem, it remains to check the limit at (\ref{c3}) holds
[simultaneously over sequences of starting points $x_n\rightarrow x$,
times $\mathbf{t}=(t_1,\ldots,t_k)$, and functions~$f$] for $\mathbb
{P}$-a.e. $\omega$. In fact, using the independent increments property
of $\widetilde{Y}^n$ and some standard analysis, it will suffice to check
the result for $k=2$ and for functions $f$ of the form $f(\widetilde
{Y}^n_{t_1},\widetilde{Y}^n_{t_2})=f_1(\widetilde
{Y}^n_{t_1})f_2(\widetilde
{Y}^n_{t_2})$. Writing the semigroup of $\widetilde{Y}^n$ as $\widetilde{P}^n$,
we have for such a function $f$ that
%
\begin{eqnarray}\label{express}
&&\biggl\llvert\int_{\overline{B}_E(x,\varepsilon)}\widetilde{E}^\omega_{ny}
f\bigl(\widetilde{Y}^n_{t_1},\widetilde{Y}^n_{t_2}
\bigr) P^\omega_{x_0} \bigl(\widetilde Y^n_{\tau^n}
\in dy|\tau^n\leq T \bigr) -\widetilde{E}^\omega_{nx_n}
f\bigl(\widetilde{Y}^n_{t_1},\widetilde{Y}^n_{t_2}
\bigr)\biggr\rrvert
\nonumber
\\
&&\qquad =\biggl\llvert\int_{\overline{B}_E(x,\varepsilon)}\widetilde{P}^n_{t_1}g(y)
P^\omega_{x_0} \bigl(\widetilde Y^n_{\tau^n}\in
dy|\tau^n\leq T \bigr) -\widetilde{P}^n_{t_1}g(x_n)
\biggr\rrvert
\\
&&\qquad \leq\sup_{y\in\overline{B}_E(x,\varepsilon)\cap n^{-1}\widetilde
{\mathcal
{C}}_1}\bigl\llvert\widetilde{P}^n_{t_1}g(y)-
\widetilde{P}^n_{t_1}g(x_n)\bigr\rrvert,\nonumber
\end{eqnarray}
where $g:=f_1\times(\widetilde{P}^n_{t_2-t_1}f_2)$ (which is a bounded,
continuous function).
Take $R>2$ with $x\in{B}_E(0,R/2)\cap n^{-1}\widetilde{\mathcal{C}_1}$.
For each $\lambda>1$, it holds that
\[
\widetilde{P}^n_{t}g=\widetilde{P}^{n,B_{\lambda R}}_{t}(g1_{B_{\lambda R}})
+\bigl(\widetilde{P}^{n}_{t}-\widetilde{P}^{n,B_{\lambda R}}_{t}
\bigr) (g1_{B_{\lambda
R}})+\widetilde{P}^n_{t}(g1_{B_{\lambda R}^c}),
\]
where ${B}_{s}:={B}_E(0,s)\cap n^{-1}\widetilde{\mathcal{C}_1}$. We have
\begin{eqnarray*}
&& \bigl|\bigl(\widetilde{P}^{n}_{t}-\widetilde{P}^{n,B_{\lambda R}}_{t}
\bigr) (g1_{B_{\lambda
R}}) (z)+\widetilde{P}^n_{t}(g1_{B_{\lambda R}^c})
(z)\bigr|
\\
&&\qquad \le \|g\|_\infty\widetilde{P}^{z}_{\omega}\bigl(
\tilde{\tau}^n_{B_{\lambda
R}}\leq t\bigr) +\|g\|_\infty
\widetilde{P}^{z}_{\omega}\bigl(\widetilde{Y}^n_t
\in B_{\lambda
R}^c\bigr)
\\
&&\qquad \le 2\|g\|_\infty\widetilde{P}^{z}_{\omega}\bigl(
\tilde{\tau}^n_{B_{\lambda
R}}\leq t\bigr).
\end{eqnarray*}
So, setting $\overline{B}_{x,\varepsilon}:=\overline{B}_E(x,\varepsilon
)\cap
n^{-1}\widetilde{\mathcal{C}_1}$, for large $n$ we have
\begin{eqnarray*}
&&\sup_{y\in\overline{B}_{x,\varepsilon}}\bigl\llvert\widetilde{P}^n_{t_1}g(y)-
\widetilde{P}^n_{t_1}g(x_n)\bigr\rrvert
\\
&&\qquad \le \sup_{y\in\overline{B}_{x,\varepsilon}}\bigl\llvert\widetilde
{P}^{n,B_{\lambda
R}}_{t_1}(g1_{B_{\lambda R}})
(y)-\widetilde{P}^{n,B_{\lambda
R}}_{t_1}(g1_{B_{\lambda R}})
(x_n)\bigr\rrvert
\\
&&\quad\qquad{} +4\|g\|_\infty\sup_{z\in\overline{B}_{x,2\varepsilon}}
\widetilde{P}^{z}_{\omega
}\bigl(\tilde{\tau}^n_{B_{\lambda R}}
\leq t\bigr)
\\
&&\qquad =:  I_1+4\|g\|_\infty I_2.
\end{eqnarray*}
Now let us note that Assumption~\ref{teoassumpHC} holds for $X^n$
killed on exiting $B_{\lambda R}$ when $D_n$ is replaced by
$B_{\lambda R}$
(which can be verified similar to the discussion in Section~\ref
{qipsec} for the $\mathbb{L}$ case; for this, it is useful to note that
the killing does not have any effect since points in $B_{\lambda R/2}$
are suitably far away from the boundary of $B_{\lambda R}$). Moreover,
applying the scaling relation $\tilde{q}^n_t(x,y)=n^dq^{\widetilde
{Y}}_{tn^2}(nx,ny)$, by
Proposition~\ref{teobiskpreLem} (for the full $\mathbb{Z}^d$ model)
we have,
${\llVert P^{n,B_{\lambda R}}_tg \rrVert}_{\infty,n,{\lambda
R}}\leq
c_1t^{-d/2}{\llVert g \rrVert}_{1,n,{\lambda R}}$ whenever\vspace*{1pt} we
also have $tn^2\geq
c_2(\sup_{x \in B_{\lambda R}}R_{x}\vee1\vee\sup_{x,y\in B_{\lambda
R}}2d_1(x,y))^{1/4}$,
where $R_{x}$ is defined as in Proposition~\ref{teobiskpreLem}.
Hence, a simple Borel--Cantelli argument using the tail estimate of
that proposition to control $\sup_{x \in B_{\lambda R}}R_{x}$ and Lemma
\ref{distlem} to control $\sup_{x,y\in B_{\lambda R}}2d_1(x,y)$ yields
that the above bound holds true for all large $n$, $\mathbb{P}$-a.s.
Hence, by applying Proposition~\ref{PR4}, we have $I_1\le
C_t(2\varepsilon)^{\gamma'}\|g\|_{2,n,\lambda R}$ for all $x\in
B_{R/2}$ and all large $n$, $\mathbb{P}$-a.s.
By applying the full $\mathbb{Z}^d$ version of Proposition~\ref{vague},
that is, the vague convergence of $m_n$ to a multiple of Lebesgue, it
follows that for $\mathbb{P}$-a.e. $\omega$: for $t>0$, $R>2$ and
$\lambda>1$,
\[
\lim_{\varepsilon\rightarrow0}\limsup_{n\rightarrow\infty}\sup
_{x \in
\overline B_{R/2}}I_1=0.
\]
For $I_2$, we will apply the full $\mathbb{Z}^d$ version of the exit
time bound of Proposition~\ref{teoldtqaiflfs}. In particular, this
result implies that for $\mathbb{P}$-a.e. $\omega$: for $t>0$, $R>2$
and $\lambda>1$,
\[
\lim_{\varepsilon\rightarrow0}\limsup_{n\rightarrow\infty}\sup
_{x \in
\overline{B}_{R/2}}I_2 \leq\lim_{\varepsilon\rightarrow0}\limsup
_{n\rightarrow\infty
}c_2 \Psi\bigl(c_3(\lambda
R/2-2\varepsilon)n,tn^2\bigr)= c_2e^{-c_3\lambda^2R^2/t}.
\]
(Note that a Borel--Cantelli argument that depends on the tail estimate
for $R_{x}$ of Proposition~\ref{teoldtqaiflfs} is hidden in the
inequality.) Letting $\lambda\rightarrow\infty$, this converges to 0.
Thus, we
deduce that for fixed $R>2$ and $\mathbb{P}$-a.e.-$\omega$ that: for
every sequence of starting points $x_n\in n^{-1}\widetilde{\mathcal{C}}_1$
such that $x_n\rightarrow x\in B_E(0,R/2)$, for every $0<t_1<t_2$, for
every bounded, continuous $f_1$, $f_2$, the expression at (\ref
{express}) converges to 0 as $n\rightarrow\infty$ and then
$\varepsilon
\rightarrow0$. Since there is no problem in extending this result to
allow any $x\in\mathbb{R}^d$, we have thus completed the proof of
(\ref{c3}).
\end{pf*}
\end{appendix}

\section*{Acknowledgements}
Parts of this work were completed during visits
by the first and second named authors to the Research Institute for
Mathematical Sciences, Kyoto University in early 2012 and early 2013,
respectively. Both authors kindly thank RIMS for the hospitality and
financial support.


%

\printaddresses

\end{document}